\documentclass[12pt]{article}

\usepackage{latexsym}
\usepackage{xfrac}  
\usepackage{setspace}
\usepackage{amsmath}
\usepackage{amssymb}
\usepackage{latexsym}
\usepackage{mathtools}
\usepackage{amsmath}
\usepackage{amssymb}
\usepackage{natbib}
\usepackage{wrapfig,lipsum,booktabs}
\usepackage{caption}
\usepackage{graphicx}
\usepackage{subcaption}
\usepackage{placeins}
\usepackage{caption}
\usepackage{color}
\setcitestyle{authoryear,open={(},close={)}}
\usepackage{hyperref}
\hypersetup{colorlinks,citecolor=blue}
\usepackage{graphics}
\usepackage{graphicx}
\usepackage{multirow}
\makeatletter
\def\ps@pprintTitle{%
  \let\@oddhead\@empty
  \let\@evenhead\@empty
  \let\@oddfoot\@empty
  \let\@evenfoot\@oddfoot
}
\makeatother
\makeatletter \oddsidemargin  -.5cm \evensidemargin -.5cm \textwidth
17.5cm \topmargin -.65in \textheight 24cm

\newtheorem{theorem}{Theorem}[section]

\newtheorem{lemma}[theorem]{Lemma}
\newtheorem{definition}{Definition}[section]
\newtheorem{proposition}[theorem]{Proposition}
\newtheorem{corollary}[theorem]{Corollary}
\newtheorem{remark}{Remark}[section]
\newtheorem{example}{Example}[section]

\numberwithin{equation}{section}
\numberwithin{equation}{section}

\title{On improved estimation of the larger location parameter}
\author{\small Naresh Garg$^a$\footnote
	{\baselineskip=10pt
		~garg.naresh22@gmail.com}, Lakshmi Kanta Patra$^b$  and Neeraj Misra$^a$   \\
    \small $^a$Department of Mathematics and Statistics,
    \small Indian Institute of Technology Kanpur\\
\small $^b$Department of Mathematics ,
\small Indian Institute of Technology Bhilai}
\spacing{1.15}
\begin{document}
\date{}
\maketitle
\begin{abstract}
This paper investigates the problem of estimating the larger location parameter of two general location families from a decision-theoretic perspective. In this estimation problem, we use the criteria of minimizing the risk function and the Pitman closeness under a general bowl-shaped loss function.  Inadmissibility of a general location and equivariant estimators is provided. We prove that a natural estimator (analogue of the BLEE of unordered location parameters) is inadmissible, under certain conditions on underlying densities, and propose a dominating estimator. We also derive a class of improved estimators using the Kubokawa's IERD approach and observe that the boundary estimator of this class is the Brewster-Zidek type estimator. Additionally, under the generalized Pitman criterion, we show that the natural estimator is inadmissible and obtain improved estimators. The results are implemented for different loss functions, and explicit expressions for the dominating estimators are provided. We explore the applications of these results to for exponential and normal distribution under specified loss functions. A simulation is also conducted to compare the risk performance of the proposed estimators. Finally, we present a real-life data analysis to illustrate the practical applications of the paper's findings.   \\

\end{abstract}
\section{Introduction}

The problem of estimating ordered parameters, whether the correct ordering between the parameters apriori, is known, has been extensively discussed in the literature. When the ordering among the parameters is known apriori, numerous studies have focused on estimating the smallest and the largest parameters, with contributions from \cite{blumenthal1968estimation},  \cite{MR1370413},  \cite{kumar2005james}, \cite{chang2015} \cite{patra}, \cite{chang2017estimation} and \cite{garg2021componentwise}. For a detailed review  on estimation of restricted parameter we refer to \cite{MR0326887}, \cite{MR961262} and \cite{MR2265239}. However, limited attention has been given to the estimation of the smallest and the largest parameters when the correct ordering between the parameters is unknown. This problem can be viewed as an estimation counterpart to the well-known "Ranking and Selection" problems, where a basic goal is to select the population associated with the largest (or smallest) parameter while lacking knowledge of the correct ordering among parameters (refer to \cite{dudewicz1982complete}, for an extensive bibliography on ranking and selection problems).

\vspace*{2mm}

In many real-world scenarios, such as in environmental studies, finance, or risk management, estimating the largest (or smallest) location/scale parameter is essential for assessing extreme events, outliers, or rare occurrences. For example, in environmental studies, estimating the largest (or smallest) location parameter of pollutant concentrations helps in determining critical levels at which adverse effects may occur. In finance, estimating the largest (or smallest) location parameter of stock returns allows investors to understand the potential for extreme losses or gains. In early work in this area, \cite{blumenthal1968estimation} considered estimation of the larger mean of two normal distributions having a common known variance. They proposed various estimators for the larger mean and compared their performances under the squared error loss function.  \cite{dhariyal1982estimation} extended the class of estimators proposed by \cite{blumenthal1968estimation} by introducing two new estimators for estimating the larger mean of two normal distributions. \cite{elfessi1992estimation} focused on the estimation of the smaller and larger scale parameters of two uniform distributions. They proposed improved estimators that outperformed the usual estimators based on the mean squared error criterion and the Pitman nearness criterion. \cite{misra1997estimation} consider two exponential distributions with unknown scale parameters. They dealt with estimation of the smaller and the larger scale parameters and obtained the MLEs. They showed that the MLEs are inadmissible and better estimators are derived under the squared error loss function. Most of the studies related to this problem have focused on specific distributions with independent marginals and specific loss functions. Some other contributions on the estimation of the larger and the smaller location/scale parameters can be found in \cite{kumar1993unbiased}, \cite{misra1994estimation}, \cite{mitra1994estimating}, and \cite{misra2002natural}. For a general framework, \cite{Misra2003} dealt with estimation of the largest scale parameter of $ k\; (\geq 2)$ independent and absolutely continuous scale parameter distributions (general probability scale models). Under the assumption of a monotone likelihood ratio on the probability models and the squared error loss function, they established that a natural estimator is inadmissible and obtained a dominating estimator. They also provided applications of these results to some specific probability models.  In this paper, we make an attempt to unify/extend various studies by considering estimation of the larger location parameters of two general probability models under a general loss function.\\



In most of the above studies, criterion of minimizing the risk function is used to obtain estimators outperforming usual estimators (such as those based on the component-wise best location/scale equivariant or the maximum likelihood estimators) under the squared error loss function. A popular alternative criterion to compare different estimators is the Pitman nearness (PN) criterion, due to \cite{pitman1937}. It compares two estimators based on the probability of one estimator being closer to the estimand than the other estimator under the absolute error loss function. \cite{rao1981} has pointed out some advantages of the Pitman nearness (PN) criterion over the mean squared error criterion. Keating (\citeyear{keating1985}) further supported Rao's observations through certain estimation problems, and \cite{keating1985m} provided some practical examples where the PN criterion is more relevant than minimizing the risk function. Additionally, \cite{peddada1985} and \cite{MR860477} extended the PN criterion to the generalized Pitman Criterion (GPN) by considering the general loss function instead of the absolute error loss function. A detailed description of the PN criterion  and the relevant literature can be found in the monograph by \cite{keating1993}. \vspace*{2mm}

The PN criterion has been extensively used in the literature for comparing various estimators in different estimation problems. However, there are only limited number in studies on the use of the PN criterion following Stein (1964) approach to obtain improvements over the usual estimators (\cite{nayak1990} and \cite{kubokawa1991}). Moreover, all these studies are centred around specific probability distributions (mostly, normal and gamma) and absolute error loss in the PN criterion. In this paper, we consider the problem of estimation of the larger location parameter of two general location models under the GPN criterion. We develop a result that is useful in finding improvements over location equivariant estimators in certain situations. 

Throughout, $\Re$ will denote the real line and $\Re^2=\Re \times \Re$ will denote the two dimensional Euclidean space. for any two real numbers $x_1$ and $x_2$, $x_{(1)}=\min\{x_1,x_2\}$ and $x_{(2)}=\max\{x_1,x_2\}$ also denote the smaller and larger of them respectively. 
\vspace*{2mm}

Let $X_1$ and $X_2$ be two independently distributed random variable with densities $f(x_1-\theta_1),\; x_1\in \Re,\; \theta_1\in \Re,$ and $f(x_2-\theta_2),\; x_2\in \Re,\; \theta_2\in \Re,$ respectively. Our aim is to estimate the larger location parameter $\theta_{(2)}=\max\{\theta_1,\theta_2\}$  under a non-negative loss function $L((\theta_1,\theta_2),a)$, $(\theta_1,\theta_2) \in \Theta=\Re^2$ and $a\in \mathcal{A}=\Re$, here $\Theta$ and $\mathcal{A}$ denotes the parametric space and the action space, respectively. At first we will invoke the principle of invariance under a suitable group of transformation. For this purpose, we consider the group $\mathcal{G}$ of transformation, where 
\begin{eqnarray*}
\mathcal{G}=\{g_c:c\in \mathbb{R}\}\cup \{g^*_1,g^*_2\},
\end{eqnarray*}
\begin{eqnarray*}
g_c(x_1,x_2)=(x_1+c,x_2+c), ~~g^*_1(x_1,x_2)=(x_1,x_2)~\mbox{ and }~g^*_2(x_1,x_2)=(x_2,x_1),\;\; x_1\in \Re,\; x_2\in \Re,\; c\in \Re.
\end{eqnarray*}
It can be easily seen that, under the group of transformation the family of distributions under consideration is  invariant and induced group of transformation of the parametric space $\Theta$ and the action space $\mathcal{A}$ as 
$$\bar{\mathcal{G}}=\{\bar{g}_c~:~ c\in \Re\}~~\mbox{ and }~~\tilde{\mathcal{G}}=\{\tilde{g}_c~:~c\in \Re\}$$
respectively, where for every $(x_1,x_2) \in \Theta$ and $c\in \Re$
\begin{eqnarray*}
\overline{g}_c(x_1,x_2)=(x_1+c,x_2+c), ~~\overline{g}^*_1(x_1,x_2)=(x_1,x_2)~\mbox{ and }~\overline{g}^*_2(x_1,x_2)=(x_2,x_1)
\end{eqnarray*}
and for any $a \in \mathcal{A}$ and $c \in \Re$ 
\begin{eqnarray*}
\widetilde{g}_c(a)=a+c, ~~\widetilde{g}_1^*(a)=a~~\mbox{ and }~~\widetilde{g}_2^*(a)=a
\end{eqnarray*}
Now the loss function $L((\theta_1,\theta_2),a)$ is invariant under the group $\mathcal{G}$ if and only is for any $(x_1,x_2) \in \Theta$, $a \in \mathcal{A}$ and $c\in \Re$
\begin{eqnarray}\label{inv1}
L(\overline{g}_c(x_1,x_2),\widetilde{g}_c(a))=L((x_1,x_2),a), ~\mbox{that is}~
L((x_1+c,x_2+c),a+c)=L((x_1,x_2),a),
\end{eqnarray}
and
\begin{eqnarray}\label{inv2}
L(\overline{g}^*_2(x_1,x_2),\widetilde{g}_2^*(a))=L((x_1,x_2),a), ~\mbox{that is}~
L((x_2,x_1),a)=L((x_1,x_2),a).
\end{eqnarray}
So combining conditions (\ref{inv1}) and (\ref{inv2}) we have, for any function $V:\Re^2\rightarrow [0,\infty)$ and for all $(x_1, x_2)\in \Theta$, $a \in \mathcal{A}$ and $c\in \Re$
\begin{eqnarray*}
L((x_1,x_2),a)=V((x_{(1)},x_{(2)}),a)~\mbox{and}~ V((x_{(1)},x_{(2)}),a)=
V((x_{(1)}+c,x_{(2)}+c),a+c).
\end{eqnarray*}
This suggests us to consider the loss function as
\begin{eqnarray}\label{loss1}
L((\theta_1,\theta_2),a)=W(a-\theta_{(2)}),\; \;(\theta_1,\theta_2)\in \Theta,\; a\in \Re,
\end{eqnarray}
where $W:\Re\rightarrow [0,\infty)$ is given function. 
Now onwards we denote $\boldsymbol{\theta}=(\theta_1,\theta_2)$.  Further we will make the following assumptions on the function $W(.)$: 
\begin{itemize}
\item	[\textbf{\boldmath(C1):}] $W(0)=0$, $W(t)$ is strictly decreasing on $(-\infty,0)$ and is strictly increasing on $(0,\infty)$, that is, $W(t)$ is strictly bowl shaped function in $t\in (-\infty,\infty)$; 
\item	[\textbf{\boldmath(C2):}] $W'(t)$ is increasing, almost everywhere;
\item 	[\textbf{\boldmath(C3):}] Integrals involving $W(t)$ are finite and differentiation under the integral sign is permissible.
\end{itemize}
\vspace*{2mm}

\noindent An estimator $\delta(x_1,x_2)$ is invariant under the group of transformation $\mathcal{G}$ if, and only if for any $(x_1,x_2)\in \Re$ and $c\in \Re$
\begin{eqnarray}\label{inavriance}
\delta(x_1+c,x_2+c)=\delta(x_1,x_2)+c~~\mbox{ and }~~\delta(x_1,x_2)=\delta(x_2,x_1).
\end{eqnarray}
From the second condition of (\ref{inavriance}), we have $\delta(x_1,x_2)=\delta^*(x_{(1)},x_{(2)})$, for some function $\delta^*$, and from the first condition of (\ref{inavriance}) we have $\delta^*(x_{(1)}+c,x_{(2)}+c)=\delta^*(x_{(1)},x_{(2)})+c$. This suggests that any invariant estimator has the form $\delta_{\phi}(X_1,X_2)$, where 
 \begin{equation}\label{investi}
 	\delta_{\phi}(X_{(1)},X_{(2)})=X_{(2)}-\phi(X_{(2)}-X_{(1)})=X_{(2)}-\phi(U)
 \end{equation}
where $U=X_{(2)}-X_{(1)}$ and $\phi:[0,\infty)\rightarrow \Re$ is real valued  function.\\

For the component problem of estimating $\theta_i$ usual estimator that is the best location equivariant estimator (BLEE) of $\theta_i$ with respect the loss function $L_i(\theta_i,\delta)=W(\delta-\theta_i)$ is obtained as
\begin{eqnarray}
	\delta^i_{c_0}(\bold{X})=X_i-c_0,\;\;i=1,2, 
\end{eqnarray}
where $c_0$ is the unique solution of equation
\begin{eqnarray}
\int_{-\infty}^{\infty}W^{\prime}(x-c_0)f(x)dx=0.
\end{eqnarray}
So, we consider a natural estimator of $\theta_{(2)}$ as
\begin{equation}\label{usual estimator}
	\delta_{c_0}(\bold{X})=X_{(2)}-c_0.
\end{equation}
Our aim is to find estimators which improve upon the natural estimator $\delta_{c_0}$, for estimating $\theta_{(2)}$ under the loss function $L((\theta_1,\theta_2),a)$ defined by (\ref{loss1}).\\

The rest of the paper is organized as follows. Inadmissibility of the usual estimator $\delta_{c_0}$ has been proved in Section \ref{sec2}. We have obtained \cite{stein1964}-type dominating estimator to demonstrate the inadmissibility. 
In Section \ref{sec3}, we consider a class of natural estimators for estimating $\theta_{(2)}$ as $\mathcal{D}=\{\delta_{b}=X_{(2)}-b: \;b\in \Re\}$. We obtain admissible estimators within the class $\mathcal{D}$. It is seen that under the condition $\lim_{t\rightarrow \infty} f(t)=0$, one of the boundary estimators of this admissible class is $\delta_{c_0}$. Furthermore, we derive a class of improved estimators over a boundary estimator of the class of admissible estimators using the IERD approach proposed by \cite{kubokawa1994unified}. Additionally, we obtain a \cite{brewster1974improving}-type improved estimator, which improves upon the boundary estimator of the admissible class of estimators. In Section \ref{sec4}, we have obtained improved estimators for the special loss function: the squared error loss, the linex loss, and the absolute error loss.\\

In section \ref{sec5}, we consider the estimation of $\theta_{(2)}$ with respect to the Pitman closeness criterion. Under the Pitman closeness, we obtain the \cite{stein1964}-type estimator, which dominates the usual estimator of $\theta_{(2)}$. As an application, in section \ref{sec6}, we derive improved estimators under the squared error loss, the linex loss, and the absolute error loss functions for the normal distribution and exponential distribution. We observe that the \cite{stein1964}-type improved estimator, and the usual estimator for the normal distribution with respect to the squared error loss and the absolute error loss are the same due to the symmetric nature of these loss functions. In Section \ref{sec7}, a simulation is carried out to compare the risk performance of the proposed estimators. Section \ref{sec8} presents a real-life data analysis, showcasing the practical applications of our paper's findings. Lastly, in Section \ref{sec9}, we offer our concluding remarks on the paper's contributions.

\section{Inadmissibility of usual estimator}\label{sec2}
In this section, we will prove inadmissibility of the natural estimator $\delta_{c_0}=X_{(2)}-c_0$ by deriving dominating estimators under the bowl shaped $\mathcal{G}$-invariant loss function \eqref{loss1}, where $W(\cdot)$ satisfies assumptions (C1)-(C3).  
\noindent To prove the dominance result we require the following assumptions:
\begin{itemize}
	\item [(A)] The family $\{f(x-\eta):\eta\in \Re\}$ of p.d.f.s holds the MLR property, i.e., for $-\infty<x_1<x_2<\infty$ and $-\infty<\eta_1<\eta_2<\infty$,  $f(x_1-\eta_1)f(x_2-\eta_2)\geq f(x_1-\eta_2)f(x_2-\eta_1)$.
\end{itemize}
\begin{itemize}
	\item [(B)] For every fixed $u>0$ and $\theta\geq 0$, let $c\equiv c(\theta,u)$ be the unique minimizer of the following function
	\begin{equation*}
		\frac{\int_{-\infty}^{\infty}W(z-c)[f(z-u+\theta)f(z)+f(z-u)f(z+\theta)]dz}
		{\int_{-\infty}^{\infty}[f(z-u+\theta)f(z)+f(z-u)f(z+\theta)]dz}.
	\end{equation*}
	That implies, for every fixed $u>0$ and $\theta\geq 0$, $c\equiv c(\theta,u)$ is the unique solution of the equation
	\begin{eqnarray*}
	\int_{-\infty}^{\infty}W^{'}(z-c(\theta,u))[f(z-u+\theta)f(z)+f(z-u)f(z+\theta)]dz
		=0.
	\end{eqnarray*}
	
\end{itemize}
\begin{theorem}\label{stth1}
	Let assumptions (A), (B), and (C1)-(C3) hold. For any fixed $u>0$, let 
	$c \equiv c(0,u)$ be the unique solution of the equation
	\begin{eqnarray}\label{eq: 2.1}
		\int_{-\infty}^{\infty}W^{'}(z-c)f(z-u)f(z)dz=0.
	\end{eqnarray}
	Then the estimator,
	\begin{eqnarray}\label{deltast1}
		\delta_{\phi_0}(\bold{X})=
			X_{(2)}-\min\{\phi(U),c(0,U)\}		
	\end{eqnarray}
	improves upon the equivariant estimator $\delta_{\phi}(\bold{X})=X_{(2)}-\phi(U)$ under the loss (\ref{loss1}), provided $P_{\boldsymbol{\theta}}\left(\phi(U) > c(0,U)\right) >0, $ at least for some $\boldsymbol{\theta}\in \Theta$.	
\end{theorem}

\noindent \textit{\textbf{Proof:}} The risk function of any equivariant estimator $\delta_{\phi}(\bold{X})=X_{(2)}-\phi(U)$ is
\begin{eqnarray*}
	R(\boldsymbol{\theta},\delta_{\phi})&=&E_{\boldsymbol{\theta}}\left[W(X_{(2)}-\phi(U)-\theta_{(2)})\right]\\
	&=&E_{\boldsymbol{\theta}}^U\left[E_{\boldsymbol{\theta}}^{X/U}\left[W\left(X_{(2)}-\phi(U)-\theta_{(2)}\right)|U\right]\right], \;\;\boldsymbol{\theta}\in \Theta.
\end{eqnarray*}
For any fixed $\boldsymbol{\theta}\in \Theta$ and $u>0$, consider
\begin{eqnarray*}
	R_{\boldsymbol{\theta},u}(c)&=&E_{\boldsymbol{\theta}}^{X/U}\left[W\left(X_{(2)}-c-\theta_{(2)}\right)|U=u\right]\\
	&=& \frac{\int_{-\infty}^{\infty}W(y-\theta_{(2)}-c)[f(y-u-\theta_{(1)})f(y-\theta_{(2)})
		+f(y-u-\theta_{(2)})f(y-\theta_{(1)})]dy}{\int_{-\infty}^{\infty}[f(y-u-\theta_{(1)})f(y-\theta_{(2)})
		+f(y-u-\theta_{(2)})f(y-\theta_{(1)})]dy}\\
		&=& \frac{\int_{-\infty}^{\infty}W(z-c)[f(z-u+\theta)f(z)+f(z-u)f(z+\theta)]dz}
		{\int_{-\infty}^{\infty}[f(z-u+\theta)f(z)+f(z-u)f(z+\theta)]dz},\;\; -\infty<c<\infty,
\end{eqnarray*}
here $\theta=\theta_{(2)}-\theta_{(1)}\in [0,\infty)$.
 	Using the assumption (B), for every fixed $\theta \geq 0$ and $u>0$, there exists a unique minimizer of $R_{\boldsymbol{\theta},u}(c)$ say $c\equiv c(\theta,u)$ which is the unique solution of the equation $R^{\prime}_{\boldsymbol{\theta},u}(c)=0$, i.e.,
\begin{eqnarray*}
\frac{\int_{-\infty}^{\infty}W^{'}(z-c(\theta,u))[f(z-u+\theta)f(z)+f(z-u)f(z+\theta)]dz}
	{\int_{-\infty}^{\infty}[f(z-u+\theta)f(z)+f(z-u)f(z+\theta)]dz}=0.
\end{eqnarray*}
For any fixed $\theta\ge 0$ and $u>0$, let $Z_{\theta,u}$ be a random variable having the density $$\Pi_{\theta,u}(z) = \frac{f(z-u+\theta)f(z)+f(z-u)f(z+\theta)}{\int_{-\infty}^{\infty}f(t-u+\theta)f(t)+f(t-u)f(t+\theta)dt},\;-\infty<z<\infty,$$
so that $E[W^{'}(Z_{\theta,u}-c(\theta,u))]=0.$ Then, for any fixed $\theta\geq 0$ and $u>0$,
\begin{eqnarray*}
	\frac{\Pi_{\theta,u}(z)}{\Pi_{0,u}(z)}&=&d(\theta,u)\frac{f(z-u+\theta)f(z)+f(z-u)f(z+\theta)}
	{2f(z-u)f(z)}\\
	&=&\frac{1}{2}\left[\frac{f(z-u+\theta)}{f(z-u)}+\frac{f(z+\theta)}{f(z)}\right],\; -\infty<z<\infty,
\end{eqnarray*}
where $d(\theta,u)$ is a positive constant.
\\By the assumption (A), we have $\displaystyle\frac{\Pi_{\theta,u}(z)}{\Pi_{0,u}(z)}$ decreasing in $z$, for any fixed $\theta\geq 0$ and $u>0$. Since, for any constant $c$, $W^{'}(z-c)$ is an almost everywhere increasing function of $z$, we conclude that, for any $u>0$,
\begin{equation}\label{eq: 2.3}
E[W^{'}(Z_{\theta,u}-c)] \le E[W^{'}(Z_{0,u}-c)],\;\forall \; \theta\geq 0,\; c\in \Re.
\end{equation}
Taking $c=c(\theta,u)$ in \eqref{eq: 2.3}, we have, for any $\theta \geq 0$ and $u>0$,
\begin{align*}
	0=&E[W^{'}(Z_{\theta,u}-c(\theta,u))] \le E[W^{'}(Z_{0,u}-c(\theta,u))]\\
	\implies\qquad  0=&E[W^{'}(Z_{0,u}-c(0,u))] \le E[W^{'}(Z_{0,u}-c(\theta,u))]
\end{align*}
Since, for any fixed $t$, $W^{'}(t-c)$ is a decreasing function of $c\in \Re$, we get
\begin{eqnarray*}
	c(\theta,u)\le c(0,u),\;\;\forall\;u>0,\; \theta\geq 0.
\end{eqnarray*}

\vspace*{2mm}

Now consider the function $\phi_0(u)=\min\{\phi(u),c(0,u)\},\; u>0$. Then, for any fixed $\theta \geq 0$ and $u>0$, we have $ c(\theta,u)\le \phi_0(u) < \phi(u)$, provided $\phi(u)>c(0,u)$. 
Using condition (C1), for any fixed $\theta \geq 0$ and $u>0$, $R_{\boldsymbol{\theta},u}(c)$ is increasing in $c\in [c(0,u),\infty)$. Consequently we get
\begin{eqnarray*}
	E_{\boldsymbol{\theta}}^{X/U}\left[W\left(\delta_{\phi_0}-\theta_{(2)}\right)|U=u\right]\le E_{\boldsymbol{\theta}}^{X/U}\left[W\left(\delta_{\phi}-\theta_{(2)}\right)|U=u\right]
\end{eqnarray*}
for all $\boldsymbol{\theta}\in\Theta$ and $u>0$ and strict inequity holds for some $u>0$. Hence we have $R(\boldsymbol{\theta},\delta_{\phi_0}) \le R(\boldsymbol{\theta},\delta_{\phi})$. This proves the theorem.
\hfill $\blacksquare$

\begin{corollary}
	Let the assumption (A) and assumptions (C1)-(C3) hold. Then the estimator,
	\begin{eqnarray}\label{deltast1}
		\delta_{\phi_0}(\bold{X})=
		X_{(2)}-\min\{c_0,c(0,U)\}		
	\end{eqnarray}	
improves upon the natural estimator $\delta_{c_0}(\bold{X})=X_{(2)}-c_0$ under the loss (\ref{loss1}) provided $P_{\boldsymbol{\theta}}(c_0>c(0,U))> 0$, for some $\boldsymbol{\theta}\in \Theta$.
\end{corollary}
\begin{remark}
	If $f(x)$ is decreasing in $x$ then $P_{\boldsymbol{\theta}}(c_0>c(0,U))>0$ for some $\boldsymbol{\theta}\in \Theta$.
\end{remark}
\subsection{A Class of improved estimators}\label{sec3}
\noindent

A natural class of estimators for estimating $\theta_{(2)}$ is $\mathcal{D}=\{\delta_b=X_{(2)}-b: \;b\in \Re\}$. Firstly, we find admissible estimators within the class of estimators $\mathcal{D}$. The risk function of an estimator $\delta_b$ is
\begin{align} 
\nonumber R(\boldsymbol{\theta},\delta_b)	&=E_{\boldsymbol{\theta}}[W(X_{(2)}-\theta_{(2)}-b)]\\ \nonumber
	&=\iint \limits _{-\infty<x_1\leq x_2<\infty} W(x_2-\theta_{(2)}-b)\left[f(x_1-\theta_{(1)})f(x_2-\theta_{(2)})+f(x_1-\theta_{(2)})f(x_2-\theta_{(1)})\right] dx_1\,dx_2 \\ \nonumber
	&=\int_{-\infty}^{\infty}\int_{-\infty}^{z}W(z-b) \left[f(x+\theta) f(z)+f(x)f(z+\theta)\right] dx\,dz\\
	&=E_{\theta}[W(Z-b)],\;\;\theta\geq 0,   \label{risk function}
\end{align}
where $Z$ is a r.v. with the density $g_{\theta}(z)=F(z+\theta)f(z)+F(z)f(z+\theta),\;z\in \Re,\; \theta\geq 0$, and $F(z)=\int_{-\infty}^{z} f(t)\,dt,\;z\in \Re$. Using the assumption (A), it is easy to verify that, for every $\theta\geq 0$, $g_{\theta}(z)/g_{0}(z)$ is decreasing in $z$.

\vspace*{2mm}
\noindent Let $b_{\theta}$ be the continues function and be the unique solution of the equation $E_{\theta}[W^{'}(Z-b)]=\int_{-\infty}^{\infty} W^{'}(x-b) g_{\theta}(x) dx=0$. Since, for every $\theta\geq 0$, $g_{\theta}(x)/g_{0}(x)$ is decreasing in $x$ and $W^{'}(x-b)$ is decreasing in $b$,  and under the assumption $\lim_{t\to \infty}f(t)=0$, it can easy to see that 
$$\inf_{\theta \geq 0} b_{\theta}=b_{\infty}=c_0\,\leq \,b_{\theta}\,\leq \, b_{0},\; \; \forall\; \theta\geq 0.$$
\begin{theorem} Suppose that the assumption (A) holds and $\lim_{t\to \infty}f(t)=0$. Then the estimators that are admissible within the class $\mathcal{D}$ are $\{X_{(2)}-b:\;  b_{\infty}\,\leq \,b_{\theta}\,\leq \, b_{0}\}$. 
\end{theorem}
\textit{\textbf{Proof:}}  Note that, for any fixed $\boldsymbol{\theta}\in\Theta$ (or fixed $\theta\geq 0$), the risk function $R(\boldsymbol{\theta},\delta)$, given by \eqref{risk function}, is uniquely minimized at $b=b_{\theta}$, it is a strictly decreasing function of $b$ on $(-\infty,b_{\theta})$ and strictly increasing function of $b$ on $(b_{\theta},\infty)$. Since, for any $\theta \geq 0$, $b_{\theta}$ is a continuous function of $\theta\in[0,\infty)$,  it assumes all values between $\inf_{\theta\geq 0} b_{\theta}=b_{\infty}=c_0$ and $\sup_{\theta\geq 0}b_{\theta}=b_{0}$, as $\theta$ varies on $[0,\infty)$. It follows that, each $b\in[b_{\infty},b_{0}]$ uniquely minimizes the risk function $R(\boldsymbol{\theta},\delta)$ at some $\boldsymbol{\theta}\in\Theta$ (or at some $\theta\geq 0$). This proves that the estimators $\{X_{(2)}-b:\; b_{\infty}\,\leq \,b_{\theta}\,\leq \, b_{0}\}$ are admissible among the estimators in the class $\mathcal{D}$.
\hfill $\blacksquare$
\vspace*{3mm}

Hence the subclass of estimators $\mathcal{D}_0=\{X_{(2)}-b:\; b_{\infty}=c_0\,\leq \,b_{\theta}\,\leq \, b_{0}\}$ is admissible within the class $\mathcal{D}$. Now, one can also consider whether improvements can be made to the estimators within the class $\mathcal{D}_0$, but it may not be possible to obtain improvements over all estimators in $\mathcal{D}_0$. Therefore, in this section, we aim to find improvements specifically for the boundary estimator $\delta_{b_0}(\bold{X}) = X_{(2)}-b_0$ within the class $\mathcal{D}_0$, where $b_0$ is the unique solution of equation
\begin{eqnarray}\label{eq:2.5}
	\int_{-\infty}^{\infty}W^{\prime}(x-b_0)f(x)\,F(x)\,dx=0,
\end{eqnarray}
where $F(x)=\int_{-\infty}^{x} f(y)\,dy,\;x\in \Re$.
\vspace*{2mm}

Now, we use the IERD approach of Kubokawa (1994) to propose a class of estimators dominating over the estimator $\delta_{b_0}$. Further, we will obtain the Brewster-Zidek (1974) type estimator improving over $\delta_{b_0}$. Consider estimation of $\theta_{(2)}$ under the loss function (1.2). Assume that the function $W(\cdot)$ is absolute continuous and satisfies the assumptions (C1), (C2) and (C3). 
\vspace*{2mm}

The following two lemmas will be useful in proving  the next result. The lemma stated in the following lemma follows from relationship between the likelihood ratio order and the revised failure rate order in the theory of stochastic orders (see \cite{MR2265633}). 


\vspace*{2mm}

\noindent The proof of the following lemma, being straightforward, is also omitted.
\begin{lemma}\label{garglemma}
	Let $s_0\in \Re$ and let $M:\Re\rightarrow\Re$ be such that $M(s)\leq 0,\; \forall \; s<s_0, $ and $M(s)\geq 0,\; \forall \; s> s_0$. Let $ M_i:\Re\rightarrow [0,\infty), \; i=1,2,$ be non-negative functions such that
	$M_1(s) M_2(s_0) \geq (\leq)\, M_1(s_0) M_2(s),\; \forall \; s<s_0,
	\text{ and } M_1(s) M_2(s_0) \leq\,(\geq)\; M_1(s_0) M_2(s),\; \forall \; s$ $>s_0.$
	Then, 
	$$ M_2(s_0) \int\limits_{-\infty}^{\infty} M(s) \, M_1(s) ds\leq\;(\geq)\; M_1(s_0) \int\limits_{-\infty}^{\infty} M(s) \, M_2(s) ds.$$
\end{lemma}

\noindent	In the following theorem, we provide a class of estimators that improve upon the natural estimator $\delta_{b_0}$.
\begin{theorem}\label{kubothm}
	Suppose that the assumption (A) holds. Additionally, assume that $W(\cdot)$ is absolutely continuous and satisfies (C1), (C2) and (C3). Let $\delta_{\phi}(\bold{X})=X_{(2)} - \phi(U)$ be a location equivariant estimator of $\theta_{(2)}$ such that 
	\begin{itemize}
		\item [(i)] $\phi(t)$ is increasing in $t\in [0,\infty)$, 
	\item[(ii)]$\lim_{t\to\infty} \phi(t)=b_0$  \item [(iii)]$\int_{-\infty}^{\infty} W^{'}(z-\phi(t))\;[F(z)-F(z-t)]f(z)\,dz \, \leq \,0,\; \forall\; t\in [0,\infty).$
	\end{itemize} Then, 
the estimator $\delta_{\phi}(\bold{X})=X_{(2)} - \phi(U)$ is an improvement
over the estimator $\delta_{b_0}(\bold{X})=X_{(2)}-b_0$.
\end{theorem}

\noindent \textbf{Proof:} Let us fix $\boldsymbol{\theta}\in\Theta$ and let $\theta=\theta_{(2)}-\theta_{(1)}$, so that $\theta \geq 0$. Consider the risk difference
\\~\\ $\Delta(\boldsymbol{\theta})$
\begin{align*}
	&= R(\boldsymbol{\theta},\delta_{b_0})-R(\boldsymbol{\theta},\delta_{\phi})\\
	& = E_{\boldsymbol{\theta}}[W(X_{(2)}-\theta_{(2)}-b_0)- W(X_{(2)}-\theta_{(2)}-\phi(U))] \\
	&= E_{\boldsymbol{\theta}}\left[\int_{U}^{\infty}\Big\{ \frac{d}{dt} W(X_{(2)}-\theta_{(2)}-\phi(t))\Big\}\; dt\right],\\
	&= \int_{-\infty}^{\infty}\int_{u=0}^{\infty}\left[\int_{u}^{\infty}\Big\{ \frac{d}{dt} W(y-\theta_{(2)}-\phi(t))\Big\}\; dt\right][f(y-u-\theta_{(1)})f(y-\theta_{(2)})+f(y-u-\theta_{(2)})f(y-\theta_{(1)})]\,du\,dy,
\end{align*}
After changing the order of integration we have 
\\~\\ $\Delta(\boldsymbol{\theta})$
\begin{align*}
	&= \int_{t=0}^{\infty}\int_{-\infty}^{\infty}\int_{u=0}^{t}\Big\{ \frac{d}{dt} W(y-\theta_{(2)}-\phi(t))\Big\}\; [f(y-u-\theta_{(1)})f(y-\theta_{(2)})+f(y-u-\theta_{(2)})f(y-\theta_{(1)})]dudydt,\\
	&=- \!\int_{t=0}^{\infty}\phi^{\prime}(t)\!\left[\int_{-\infty}^{\infty}\int_{u=0}^{t}\!\!\Big\{  W^{\prime}(y-\theta_{(2)}-\phi(t))\Big\}\! [f(y-u-\theta_{(1)})f(y-\theta_{(2)})\!+\!f(y-u-\theta_{(2)})f(y-\theta_{(1)})]dudy\right]\!dt.
\end{align*}
Since $\phi(t)$ is a increasing function of $t$, it suffices to show that, for every $t>0$, 
\begin{align}\label{eq:2.10}	&\int_{-\infty}^{\infty}\int_{u=0}^{t}\Big\{  W^{\prime}(y-\theta_{(2)}-\phi(t))\Big\}\; [f(y-u-\theta_{(1)})f(y-\theta_{(2)})+f(y-u-\theta_{(2)})f(y-\theta_{(1)})]dudy\leq \,0   \nonumber\\
	\iff \qquad & \int_{u=0}^{t}\int_{-\infty}^{\infty}\Big\{  W^{\prime}(z-\phi(t))\Big\}\; [f(z-u+\theta)f(z)+f(z-u)f(z+\theta)]dzdu\leq \,0.
\end{align}
Now, since $W^{'}(t)$ is increasing function of $t$ and, for every fixed $\theta\geq 0$ and $u\in \Re$, $\frac{f(z-u+\theta)f(z)+f(z-u)f(z+\theta)}{2 f(z-u)f(z)}$ is decreasing in $z$, then, for $\theta\geq 0$, we have 
\begin{align*}
	&\int_{u=0}^{t}\left[\int_{-\infty}^{\infty} W^{'}(z-\phi(t))\left[f(z-u+\theta)f(z)+f(z-u)f(z+\theta)\right]dz\right]du\\
	&\leq \,\int_{u=0}^{t}\left[\int_{-\infty}^{\infty} W^{'}(z-\phi(t))\;[f(z-u)f(z)+f(z-u)f(z)]dz\right] du\\
	&= \,	\,2\int_{-\infty}^{\infty} W^{'}(z-\phi(t))\;[F(z)-F(z-t)]f(z)\,dz
\end{align*}
Now, using hypothesis (iii), we obtain \eqref{eq:2.10}. This completes proof of the theorem. 
\hfill $\blacksquare$\\

In the following we will prove a corollary which will provide \cite{brewster1974improving} type improved estimators. 

\begin{corollary}
	Suppose that assumptions (A), (C1), (C2) and (C3) hold. Additionally suppose that, for every fixed $t$, the equation
$$k_1(c\vert  t)=\int_{-\infty}^{\infty} \; W^{'}(z-c)\; [F(z)-F(z-t)] f(z)\,dz =0$$
has the unique solution $c\equiv \phi_{0}(t)$. Then
$$R(\boldsymbol{\theta},\delta_{\phi_{0}})\leq R(\boldsymbol{\theta},\delta_{0}), \;\;\; \forall \; \; \boldsymbol{\theta} \in \Theta,$$
where $\delta_{\phi_{0}}(\bold{X})=X_{(2)}-\phi_{0}(U)$.
\end{corollary}

\noindent \textbf{Proof:} It is suffices to show that $\phi_{0}(t)$ satisfies conditions of Theorem \ref{kubothm}. Note that a hypothesis of the corollary ensures that $\lim_{t\to\infty} \phi_{0}(t)=b_0$. To show that $\phi_{0}(t)$ is an increasing function of $t$, suppose that, there exist numbers $t_1$ and $t_2$ such that $0<t_1<t_2$ and $\phi_{0}(t_1)\neq \phi_{0}(t_2).$ Under the hypotheses of the corollary, we have $k_1(\phi_{0}(t_1)\vert t_1)=0$, $\phi_{0}(t_2)$ is the unique solution of $k_1(c\vert t_2)=0$ and $k_1(c\vert t_2)$ is a decreasing function of $c$. Let $ s_0 = \phi_{0}(t_1), \; M(s)=W^{'}(s-s_0)f(s),\; M_1(s)=\int_{0}^{t_2}f(s-u)du$ and $M_2(s)=\int_{0}^{t_1}f(s-u)du$. Then, using Lemma \ref{garglemma}, we get
\begin{small}
	$$\int_{0}^{t_1}f(\phi_{0,1}(t_1)-w)\,du\; \left(\int_{-\infty}^{\infty}\; W^{'}(z-\phi_{0}(t_1))\,f(z)\;\int_{0}^{t_2} f(z-u)\,du \;dz\right) \qquad \qquad \qquad \qquad \qquad \qquad \quad$$   $$\qquad \qquad \qquad \qquad \geq\; \int_{0}^{t_2}f(\phi_{0}(t_1)-w)du\; \left( \int_{-\infty}^{\infty}\; W^{'}(z-\phi_{0}(t_1))\,f(z)\;\int_{0}^{t_1} f(z-u)\,du \;dz \right)=0.$$
\end{small}
This implies that
$$k_1(\phi_{0}(t_1)\vert t_2)=\int_{-\infty}^{\infty}\int_{0}^{t_2}\; W^{'}(z-\phi_{0,1}(t_1))f(z-u)f(z)\,du\, dz\geq \; 0.$$
So we have  $k_1(\phi_{0}(t_1)\vert t_2)\,> \, 0$, as $k_1(c\vert t_2)=0$ has the unique solution $c\equiv \phi_{0}(t_2)$ and $\phi_{0}(t_1)\neq\phi_{0}(t_2)$. Since $k_1(c\vert t_2)$ is a decreasing function of $c$, $k_1(\phi_{0}(t_2)\vert t_2)$ $=0$ and $k_1(\phi_{0}(t_1)\vert t_2)\,> \, 0$, it follows that $\phi_{0}(t_1)<\, \phi_{0}(t_2)$. 
Hence the result follows. 
\hfill $\blacksquare$

\vspace*{5mm}
\section{Dominance result for special loss functions}\label{sec4}
In this section, we have obtained the improved estimators for three special loss functions namely squared error loss $L_1:W(t)=t^2,\; t\in \Re$, linex loss $L_2:W(t)=e^{at}-at-1,,\; t\in \Re,\; a\ne 0$ and absolute error loss $L_3:W(t)=|t|,\; t\in \Re$. 

\begin{theorem}\label{thml1}
	Suppose that the assumption (A) holds. Then for estimating $\theta_{(2)}$ with respect to loss function $L_1$ the estimator
	\begin{eqnarray}\label{deltastl1}
		\delta_{ST}(\bold{X})=
			X_{(2)}-\min\{c_0,c(0,u)\}
	\end{eqnarray}
	improves upon the estimator $\delta_{0}(\bold{X})$, provided $P_{\boldsymbol{\theta}}(c_0 > c(0,U)) \ne0 $, for some $\boldsymbol{\theta}\in \Theta$, where 
	\begin{eqnarray*}
		c(0,u)=\int_{-\infty}^{\infty}zf(z-u)f(z)dz\bigg/\int_{-\infty}^{\infty}f(z-u)f(z)dz.
	\end{eqnarray*}
and $c_0=\int_{-\infty}^{\infty}xf(x)dx.$

\end{theorem}

\begin{theorem}\label{thl2}
	The estimator
	\begin{eqnarray}\label{deltastl2}
		\delta_{ST}(\bold{X})=
			X_{(2)}-\min\{c_0,c(0,U)\},
	\end{eqnarray}
	improves upon the estimator $\delta_{c_0}(\bold{X})$ with respect to the loss $L_2$, where $c_0=\frac{1}{a}\ln \int_{\infty}^{\infty}e^{ax}f(x)dx$ and $c(0,U)=\frac{1}{a}\ln H(U)$
with 
\begin{eqnarray*}
H(u)=\int_{-\infty}^{\infty}e^{az}f(z-u)f(z)dz\bigg/\int_{-\infty}^{\infty}f(z-u)f(z)dz\end{eqnarray*}
provided the assumption (A) holds and $P_{\boldsymbol{\theta}}(c_0> \frac{1}{a}\ln H(u))\ne 0$, for some $\boldsymbol{\theta}\in \Theta$.
\end{theorem}

\begin{theorem}\label{thl3}
	The estimator
	\begin{eqnarray}\label{deltastl2}
		\delta_{ST}(\bold{X})=
		X_{(2)}-\min\{c_0,c(0,U)\},
	\end{eqnarray}
	improves upon the estimator $\delta_{c_0}(\bold{X})$ with respect to the loss $L_3$, where $c_0$ and $c(0,u)$ are such that
	\begin{eqnarray*}
		\int_{-\infty}^{c_0}f(z)dz=\frac{1}{2} \; \text{  and  } \;
		\int_{-\infty}^{c(0,u)}f(z-u)f(z)dz\bigg/\int_{-\infty}^{\infty}f(z-u)f(z)dz=\frac{1}{2},
	\end{eqnarray*}
	respectively, provided the assumption (A) holds and $P_{\boldsymbol{\theta}}(c_0> c(0,U))\ne 0$, for some $\boldsymbol{\theta}\in \Theta$.
\end{theorem}

\noindent Now we will apply Theorem \ref{kubothm} to particular loss functions and provide a class of improved estimators over the estimator $\delta_{b_0}(\bold{X}) = X_{(2)}-b_0,$ where $b_0$ be as defined by the equation \eqref{eq:2.5}. 

\begin{theorem}\label{kubothml1}
	Suppose that the assumption (A) holds. Let $\delta_{\phi}(\bold{X})=X_{(2)} - \phi(U)$ be a location equivariant estimator of $\theta_{(2)}$ such that 
	\begin{itemize}
		\item [(i)] $\phi(t)$ is increasing in $t$, 
		\item[(ii)]$\lim_{t\to\infty} \phi(t)=2\int_{-\infty}^{\infty}xf(x)F(x)dx=b_0,$ 
		 \item [(iii)]$\phi(t) \ge \frac{\int_{-\infty}^{\infty}\int_{0}^{t}zf(z-u)f(z)\,du\,dz}{\int_{-\infty}^{\infty}\int_{0}^{t}f(z-u)f(z)\,du\,dz}$.
	\end{itemize} Then 
	the estimator $\delta_{\phi}(\bold{X})=X_{(2)} - \phi(U)$
	improves upon the estimator $\delta_{b_0}(\bold{X})=X_{(2)}-b_0$ with respect to the $L_1$.
\end{theorem}
\begin{remark}
	The boundary estimator of the class estimators given by Theorem \ref{kubothml1} is the Brewster-Zidek type estimator. So the Brewster-Zidek type estimator is obtained as 
	$$\delta_{BZ}(\bold{X})=X_{(2)} - \phi_{BZ}(U)$$
	with 
	$$\phi_{BZ}(t)=\frac{\int_{-\infty}^{\infty}\int_{0}^{t}zf(z-u)f(z)\,du\,dz}{\int_{-\infty}^{\infty}\int_{0}^{t}f(z-u)f(z)\,du\,dz}.$$
\end{remark}

\begin{theorem}\label{kubothml2}
	Suppose that the assumption (A) holds. Let $\delta_{\phi}(\bold{X})=X_{(2)} - \phi(U)$ be a location equivariant estimator of $\theta_{(2)}$ such that 
	\begin{itemize}
		\item [(i)] $\phi(t)$ is increasing in $t$, 
		\item[(ii)]$\lim_{t\to\infty} \phi(t)=\frac{1}{a}\ln\left( 2\int_{-\infty}^{\infty}e^{ax}f(x)F(x)dx\right)=b_0,$ 
		\item [(iii)]$\phi(t) \le\frac{1}{a}\ln \left( \frac{\int_{-\infty}^{\infty}\int_{0}^{t}e^{az}f(z-u)f(z)\,du\,dz}{\int_{-\infty}^{\infty}\int_{0}^{t}f(z-u)f(z)\,du\,dz}\right)$.
	\end{itemize} Then 
	the estimator $\delta_{\phi}(\bold{X})=X_{(2)} - \phi(U)$
	improves upon the estimator $\delta_{b_0}(\bold{X})=X_{(2)}-b_0$ with respect to the $L_2$.
\end{theorem}
\begin{remark}
	The boundary estimator of the class estimators given by Theorem \ref{kubothml2} is the Brewster-Zidek type estimator.  So the Brewster-Zidek type estimator is obtained as 
	$$\delta_{BZ}(\bold{X})=X_{(2)} - \phi_{BZ}(U)$$
	with 
	$$\phi_{BZ}(t)=\frac{1}{a}\ln \left( \frac{\int_{-\infty}^{\infty}\int_{0}^{t}e^{az}f(z-u)f(z)\,du\,dz}{\int_{-\infty}^{\infty}\int_{0}^{t}f(z-u)f(z)\,du\,dz}\right),\;\;t>0.$$
\end{remark}

\begin{theorem}\label{kubothml3}
	Suppose that the assumption (A) holds and $c_0$ be as in Theorem \ref{thl3}. Let $\delta_{\phi}(\bold{X})=X_{(2)} - \phi(U)$ be a location equivariant estimator of $\theta_{(2)}$ such that 
	\begin{itemize}
		\item [(i)] $\phi(t)$ is increasing in $t$, 
		\item[(ii)]$\lim_{t\to\infty} \phi(t)=b_0$, where $b_0$ is the solution of the equation $2\int_{-\infty}^{b_0} f(x) F(x) dx=\frac{1}{2}$, 
		\item [(iii)] $\phi(t)$ be such that it satisfies
		$$\int_{-\infty}^{\phi(t)}\int_{0}^{t}f(z-u)f(z)\,du\,dz \ge \frac{1}{2} \int_{-\infty}^{\infty}\int_{0}^{t}f(z-u)f(z)\,du\,dz.$$ 
	\end{itemize} Then 
	the estimator $\delta_{\phi}(\bold{X})=X_{(2)} - \phi(U)$
	improves upon the estimator $\delta_{b_0}(\bold{X})=X_{(2)}-b_0$ with respect to the $L_3$.
\end{theorem}

\begin{remark}
	The boundary estimator of the class estimators given by Theorem \ref{kubothml3} is the Brewster-Zidek type estimator. So the Brewster-Zidek type estimator is obtained as 
	$$\delta_{BZ}(\bold{X})=X_{(2)} - C,$$
where $C$ is the unique solution of the equation
$$	\int_{-\infty}^{C}\int_{0}^{t}f(z-u)f(z)\,du\,dz = \frac{1}{2} \int_{-\infty}^{\infty}\int_{0}^{t}f(z-u)f(z)\,du\,dz,\;\;\;t>0.$$
\end{remark}

\section{Improved estimator under generalized Pitman closeness criterion}\label{sec5}
In this section we consider the estimation problem under the generalized Pitman nearness (GPN) criterion. A brief discussion on the Pitman nearness criterion is given in  \cite{garg2022}.  For completeness in our presentation, we are again discussing it here.  The notion of the Pitman nearness criterion  was first introduced by Pitman (\citeyear{pitman1937}), as defined below. 
\begin{definition}
Let $\bold{X}$ be a random vector having a probability distribution involving an unknown parameter $\boldsymbol{\theta}\in \Theta$ ($\boldsymbol{\theta}$ may be vector valued).  Let $\delta_1$ and $\delta_2$ be two estimators of a real-valued estimand $\tau(\boldsymbol{\theta})$. Then, the Pitman nearness (PN) of $\delta_1$ relative to $\delta_2$ is defined by
$$PN(\delta_1,\delta_2;\boldsymbol{\theta})=P_{\boldsymbol{\theta}}[\vert \delta_1-\tau(\boldsymbol{\theta})\vert <\vert \delta_2-\tau(\boldsymbol{\theta})\vert], \; \;\boldsymbol{\theta}\in\Theta,$$
and the estimator $\delta_1$ is said to be nearer to $\tau(\boldsymbol{\theta})$ than $\delta_2$ if $PN(\delta_1,\delta_2;\boldsymbol{\theta})\geq \frac{1}{2},\;\forall\; \boldsymbol{\theta}\in\Theta$, with strict inequality for some $\boldsymbol{\theta}\in\Theta$.
\end{definition}
 Nayak (\citeyear{nayak1990}) and Kubokawa (\citeyear{kubokawa1991}) modified the Pitman (\citeyear{pitman1937}) nearness criterion and defined the generalized Pitman nearness (GPN) criterion based on general loss function $L(\boldsymbol{\theta},\delta).$

 \begin{definition}
 	Let $\bold{X}$ be a random vector having a distribution involving an unknown parameter $\boldsymbol{\theta}\in \Theta$ and let $\tau(\boldsymbol{\theta})$ be a real-valued estimand. Let $\delta_1$ and $\delta_2$ be two estimators of the estimand $\tau(\boldsymbol{\theta})$. Also, let $L(\boldsymbol{\theta},a)$ be a specified loss function for estimating $\tau(\boldsymbol{\theta})$. Then, the generalized Pitman nearness (GPN) of $\delta_1$ relative to $\delta_2$ is defined by
$$GPN(\delta_1,\delta_2;\boldsymbol{\theta})=P_{\boldsymbol{\theta}}[L(\boldsymbol{\theta},\delta_1) <L(\boldsymbol{\theta},\delta_2)]+\frac{1}{2} P_{\boldsymbol{\theta}}[L(\boldsymbol{\theta},\delta_1) =L(\boldsymbol{\theta},\delta_2)], \; \; \boldsymbol{\theta}\in\Theta.$$
The estimator $\delta_1$ is said to be nearer to $\tau(\boldsymbol{\theta})$ than $\delta_2$, under the GPN criterion, if $GPN(\delta_1,\delta_2;\boldsymbol{\theta})\geq \frac{1}{2},\;\forall\; \boldsymbol{\theta}\in\Theta$, with strict inequality for some $\boldsymbol{\theta}\in\Theta$.
 \end{definition}

\vspace*{2mm}

The following result, popularly known as Chebyshev's inequality, will be used in our study (see \cite{MR2363282}).

\begin{proposition}
Let $S$ be random variable  and let $k_1(\cdot)$ and $k_2(\cdot)$ be real-valued monotonic functions defined on the distributional support of the r.v. $S$. If $k_1(\cdot)$ and $k_2(\cdot)$ are monotonic functions of the same (opposite) type, then
$$E[k_1(S)k_2(S)]\geq (\leq ) E[k_1(S)] E[k_2(S)],$$
provided the above expectations exist.
\end{proposition}

\vspace*{3mm}

 The following lemma, taken from \cite{garg2022}, will be useful in proving the main results of this section (also see \cite{nayak1990}) and \cite{zhou2012}).
 
\begin{lemma}[\cite{garg2022}]\label{garg1}
 Let $Y$ be a random variable having the Lebesgue probability density function and let $m_Y$ be the median of $Y$. Let $W:\Re\rightarrow [0,\infty)$ be a function such that $W(0)=0$, $W(t)$ is strictly decreasing on $(-\infty,0)$ and strictly increasing on $(0,\infty)$. Then, for $-\infty< c_1<c_2\leq m_Y$ or $-\infty<m_Y\leq c_2<c_1$,
$GPN= P[W(Y-c_2)<W(Y-c_1)]+\frac{1}{2} P[W(Y-c_2)=W(Y-c_1)]>\frac{1}{2}$.
\end{lemma}

%

%
%

%

\begin{lemma}\label{pitlemma2}
	 Suppose that assumptions (A) and (C1)-(C3) hold. For $u>0$ and $\boldsymbol{\theta}\in \Theta$, let $m(\boldsymbol{\theta},u)$ denote the median of the conditional distribution of $X_{(2)}-\theta_{(2)}$, given $U=u$. Then $$m(\boldsymbol{\theta},u)\leq m(\boldsymbol{0},u),\; \; \forall \; u>0.$$
\end{lemma}

\noindent \textit{\textbf{Proof:}}
The joint distribution of $(Y,U)=(X_{(2)}-\theta_{(2)},X_{(2)}-X_{(1)})$ is
\begin{eqnarray*}
	h_{\theta}(y,u)=f(y-u+\theta)f(y)
	+f(y-u)f(y+\theta),\;\; -\infty<y<\infty,\; u>0,
\end{eqnarray*}
where $\theta=\theta_{(2)}-\theta_{(1)}\geq 0$. The conditional p.d.f. of $Y$ given $U=u$ is $\Pi_{\theta}(z\vert u) = \frac{1}{c_{\theta}}f(z-u+\theta)f(z)+f(z-u)f(z+\theta),\,z\in \Re$, where $c_{\theta}=\int_{-\infty}^{\infty}\left[f(t-u+\theta)f(t)+f(t-u)f(t+\theta)\right] dt$. \\
 Now, for any fixed $u>0$ and $\theta\geq 0$, we have
\begin{eqnarray*}
	\frac{\Pi_{\theta}(z\vert u)}{\Pi_{0}(z\vert u)}&=&\frac{c_0}{c_{\theta}}\frac{f(z-u+\theta)f(z)+f(z-u)f(z+\theta)}
	{2f(z-u)f(z)}\\
	&=&\frac{c_0}{2\,c_{\theta}}\left[\frac{f(z-u+\theta)}{f(z-u)}+\frac{f(z+\theta)}{f(z)}\right].
\end{eqnarray*}
By the assumption (A), for $\theta>0$, $\frac{\Pi_{\theta}(z\vert u)}{\Pi_{0}(z\vert u)}$  is decreasing in $z$.

For every fixed $u>0$ and $\theta\geq 0$, take $k_1(s)=I_{(-\infty,m(\boldsymbol{\theta},u))}(s),\;s\in \Re,$ and $k_2(s)=	\frac{\Pi_{\theta}(s\vert u)}{\Pi_{0}(s\vert u)},\; \in \Re$, where $I_A(\cdot)$ denotes the indicator function of set $A  \subseteq \Re$. Here $k_1(s)$ and $k_2(s)$ are decreasing functions of $s$. Using Proposition 4.1, for any $u>$ and $\theta\geq 0$, we get
	\begin{align*}
		\frac{1}{2}=&\int_{-\infty}^{\infty} k_1(s) k_2(s) \Pi_{0}(s\vert u)ds \geq \left( \int_{-\infty}^{\infty} k_1(s) \Pi_{0}(s\vert u)ds \right) \left( \int_{-\infty}^{\infty} k_2(s) \Pi_{0}(s\vert u)ds\right)\\
		\implies &\int_{-\infty}^{m(\boldsymbol{\theta},u)} \Pi_{\theta}(s\vert u)ds=\int_{-\infty}^{m(\boldsymbol{0},u)} \Pi_{0}(s\vert u)ds=\frac{1}{2} \geq \int_{-\infty}^{m(\boldsymbol{\theta},u)} \Pi_{0}(s\vert u)ds\\
		\implies & m(\boldsymbol{0},u)\geq \, m(\boldsymbol{\theta},u),\;\; \forall \; u>0,
	\end{align*}
	establishing the assertion.
 \hfill $\blacksquare$

\begin{theorem}\label{th1}
	Suppose that assumptions (A) and (C1)-(C3) hold. Let $\delta_{\phi}=X_{(2)}-\phi(U)$ be an estimator of $\theta_{(2)}$ such that $P_{\boldsymbol{\theta}}(\phi(U)>m(\boldsymbol{0},U)) \ne 0 $, for some $\boldsymbol{\theta}\in \Theta$.
	Then the estimator
	\begin{eqnarray}\label{deltast}
		\delta_{\phi_0}(\bold{X})=X_{(2)}-\min\{\phi(U),m(\boldsymbol{0},U)\}
	\end{eqnarray}
	improves over the estimator $\delta_{\phi}(\bold{X})=X_{(2)}-\phi(U)$ in terms of the GPN criterion with a general loss ( \ref{loss1}).
	
\end{theorem}

\noindent \textit{\textbf{Proof:}}  Let $\phi_0(t)=\min\{\phi(t),m(\boldsymbol{0},t)\},\; t\geq 0.$ The GPN of $\delta_{\phi_0}(\bold{X})=X_{(2)}-\phi_0(U)$ relative to $\delta_{\phi}(\bold{X})$ is given by 
\begin{align*}
	GPN(\delta_{\phi_0},\delta_{\phi};\boldsymbol{\theta})&=P_{\boldsymbol{\theta}}[W(X_{(2)}-\theta_{(2)}-\phi_0(U)) <W(X_{(2)}-\theta_{(2)}-\phi(U))]
	\\&\quad +\frac{1}{2} P_{\boldsymbol{\theta}}[W(X_{(2)}-\theta_{(2)}-\phi_0(U)) =W(X_{(2)}-\theta_{(2)}-\phi(U))]\\
	&=\int_{-\infty}^{\infty} P_{\boldsymbol{\theta}}[W(X_{(2)}-\theta_{(2)}-\phi_0(u)) <W(X_{(2)}-\theta_{(2)}-\phi(u)) \vert U=u]\,f_U(u) \,du
	\\&\quad +\frac{1}{2}\,  \int_{-\infty}^{\infty}P_{\boldsymbol{\theta}}[W(X_{(2)}-\theta_{(2)}-\phi_0(u)) =W(X_{(2)}-\theta_{(2)}-\phi(u)) \vert U=u] \,f_U(u) \,du,
\end{align*}
where $f_U(u)=\int_{-\infty}^{\infty}\left[f(t-u+\theta)f(t)+f(t-u)f(t+\theta)\right] dt,\;u>0,\;\theta=\theta_{(2)}-\theta_{(1)}\geq 0,$ is p.d.f. of r.v. $U$. Now, define $A=\{u>0: \phi(u)\leq m(\bold{0},u)\}$ and $B=\{u>0: \phi(u)> m(\bold{0},u)\}$, we get
\begin{align*}
	GPN(\delta_{\phi_0},\delta_{\phi};\boldsymbol{\theta})&=\int_{A} P_{\boldsymbol{\theta}}[W(X_{(2)}-\theta_{(2)}-\phi(u)) <W(X_{(2)}-\theta_{(2)}-\phi(u)) \vert U=u]\,f_U(u) \,du
	\\&\quad +\frac{1}{2}\,  \int_{A}P_{\boldsymbol{\theta}}[W(X_{(2)}-\theta_{(2)}-\phi(u)) =W(X_{(2)}-\theta_{(2)}-\phi(u)) \vert U=u] \,f_U(u) \,du
	\\ &\quad + \int_{B} P_{\boldsymbol{\theta}}[W(X_{(2)}-\theta_{(2)}-m(\bold{0},u)) <W(X_{(2)}-\theta_{(2)}-\phi(u)) \vert U=u]\,f_U(u) \,du
	\\&\quad +\frac{1}{2}\,  \int_{B}P_{\boldsymbol{\theta}}[W(X_{(2)}-\theta_{(2)}-m(\bold{0},u)) =W(X_{(2)}-\theta_{(2)}-\phi(u)) \vert U=u] \,f_U(u) \,du\\
	&=\frac{1}{2}\,  \int_{A}\,f_U(u) \,du
	\\ &\quad + \int_{B} P_{\boldsymbol{\theta}}[W(X_{(2)}-\theta_{(2)}-m(\bold{0},u)) <W(X_{(2)}-\theta_{(2)}-\phi(u)) \vert U=u]\,f_U(u) \,du.
\end{align*}
\noindent  From Lemma \ref{pitlemma2}, we have
\begin{eqnarray*}
	m(\boldsymbol{\theta},u)\le m(\boldsymbol{0},u),\;\;\forall\;u>0.
\end{eqnarray*}
Then, for every $\boldsymbol{\theta}\in \Theta$, we have $-\infty<m(\boldsymbol{\theta},u)\le m(\boldsymbol{0},u) < \phi(u)$,$\; \forall\; u\in B$. Since, for $u\in B$, $m(\boldsymbol{\theta},u)$ is the median of the conditional distribution of $X_{(2)}-\theta_{(2)}$ given $U=u$ and using Lemma \ref{garg1}, we have $P_{\boldsymbol{\theta}}[W(X_{(2)}-\theta_{(2)}-m(\bold{0},u)) <W(X_{(2)}-\theta_{(2)}-\phi(u)) \vert U=u]\geq \frac{1}{2},\;\forall\; u\in B, \; \boldsymbol{\theta}\in \Theta.$ Hence we get
\begin{align*}GPN(\delta_{\phi_0},\delta_{\phi};\boldsymbol{\theta})\,
	\geq \,\frac{1}{2},\;\forall \; \boldsymbol{\theta}\in\Theta, \;u>0,
\end{align*}
and strict inequity holds for some $u>0$. This proves the theorem.
\hfill $\blacksquare$

\vspace*{3mm}

An immediate consequence of Lemma \ref{garg1}, the Pitman nearest (PN) equivariant estimator of $\theta_{(2)}$ within the class $\mathcal{D}$, under the GPN criterion, is obtained as $$\delta_{PN}(\bold{X})=X_{(2)}-m_{0},$$
where $m_{0}$ is such that $\int_{-\infty}^{m_0} f(x)\,dx=\frac{1}{2}$.

\begin{corollary}\label{pitman corr}
	Suppose that assumptions (A) and (C1)-(C3) hold. Then for estimating $\theta_{(2)}$ under the GPN criterion, the estimator
	\begin{eqnarray}\label{deltast1}
		\delta^*_{PN}(\bold{X})=X_{(2)}-\min\{m_0,m(\boldsymbol{0},U)\}
	\end{eqnarray}
	is Pitman nearer to $\theta_{(2)}$ than the estimator $\delta_{PN}(\bold{X})=X_{(2)}-m_{0}$, provided $P_{\boldsymbol{\theta}}(m(\boldsymbol{0},U)<m_0) \ne0 $, for some $ \boldsymbol{\theta}\in \Theta$.
\end{corollary}
The following corollary provides an improvement over the usual estimator $\delta_{0}(\bold{X})=X_{(2)}-c_{0}$ (as defined by \eqref{usual estimator}).
\begin{corollary}\label{pitman corr over usual estimator}
	Suppose that assumptions (A) and (C1)-(C3) hold. Then for estimating $\theta_{(2)}$ under the GPN criterion, the estimator
	\begin{eqnarray}\label{deltast1}
		\delta^*_{PN}(\bold{X})=X_{(2)}-\min\{c_0,m(\boldsymbol{0},U)\}
	\end{eqnarray}
	is Pitman nearer to $\theta_{(2)}$ than the estimator $\delta_{0}(\bold{X})=X_{(2)}-c_{0}$, provided $P_{\boldsymbol{\theta}}(m(\boldsymbol{0},U)<c_0) \ne 0 $, for some $ \boldsymbol{\theta}\in \Theta$.
\end{corollary}
\section{Applications}\label{sec6}
\begin{example} \rm
	Let $X_1\sim N(\theta_1,\sigma^2)$ and $X_2\sim N(\theta_2,\sigma^2)$, where $(\theta_1,\theta_2)\in \Theta_0$ is vector of unknown means and $\sigma^2>0$ is known common variance. The pdf of $X_i$ is $f(x-\theta_i)=\frac{1}{\sqrt{2 \pi} \sigma} e^{-\frac{1}{2\sigma^2}(x-\theta_i)^2},\;x\in \Re,\;i=1,2.$ Here it is easy to see that $f(\cdot)$ holds the assumption (A).
\vspace*{2mm}

\noindent Consider the estimation of parameter $\theta_{(2)}=\max\{\theta_1,\theta_2\}$ under the loss function
\begin{eqnarray}\label{eq:4.1}
	L((\theta_1,\theta_2),a)=W(a-\theta_{(2)}), \; (\theta_1,\theta_2)\in \Re^2,\; a\in \Re.
\end{eqnarray}
For the squared error loss (i.e., $W(t)=t^2,\;t\in \Re$), the usual estimator of $\theta_{(2)}$ is $\delta_{c_0}=X_{(2)}$. In this case, using Theorem 3.1, there is no improvement over the usual estimator $\delta_{c_0}$ and the improved estimator $\delta_{ST}(\bold{X})=X_{(2)}-\min\bigg\{0,\frac{U}{2}\bigg\}=X_{(2)}$ is same. Also, using Theorem 3.4, the estimator $\delta_{b_0}(\bold{X})=X_{(2)}-\frac{\sigma}{\sqrt{\pi}}$ is dominated by the estimator $\delta_{BZ}(\bold{X})$, where
\begin{align*}
		\delta_{BZ}(\bold{X})&=X_{(2)}-\phi_{BZ}(U)=X_{(2)}-\frac{\sigma}{\sqrt{2}} \frac{\left(\frac{1}{\sqrt{2\pi}}-\phi\left(\frac{U}{\sqrt{2} \sigma}\right)\right)}{\left(\Phi\left(\frac{U}{\sqrt{2} \sigma}\right)-0.5\right)},
\end{align*}
$U=X_{(2)}-X_{(1)}$, $\phi(\cdot)$ is the p.d.f. of standard normal distribution and $\Phi(\cdot)$ is the d.f. of standard normal distribution.
 
\vspace*{2mm}
For $W(t)=e^{at}-at-1,\;t\in \Re,\; a\neq 0$, (i.e., the Linex loss) the usual estimator of $\theta_{(2)}$ is $\delta_{c_0}=X_{(2)}-\frac{a \sigma^2}{2}$. In this case, using Theorem 3.2, the estimator $\delta_{c_0}$ is dominated by the estimator $\delta_{ST}(\bold{X})$, where
\begin{align*}
	\delta_{ST}(\bold{X})&=X_{(2)}-\min\bigg\{\frac{a \sigma^2}{2},\frac{U}{2}+\frac{a \sigma^2}{4}\bigg\}=\max\bigg\{X_{(2)}-\frac{a \sigma^2}{2},\frac{X_1+X_2}{2}-\frac{a \sigma^2}{4}\bigg\}.
\end{align*}
Also, using Theorem 3.5, the estimator $\delta_{b_0}(\bold{X})=X_{(2)}-\frac{1}{a}\left[\ln2+\frac{a^2\sigma^2}{2} +\ln\left(\Phi\left(\frac{a\sigma}{\sqrt{2}}\right)\right)\right]$ is dominated by the estimator $\delta_{BZ}(\bold{X})$, where
$$	\delta_{BZ}(\bold{X})=X_{(2)}-\frac{1}{a} \left[\frac{a^2\sigma^2}{2} +\ln\left(\Phi\left(\frac{a\sigma}{\sqrt{2}}\right)-\Phi\left(\frac{-U+a\sigma^2}{\sqrt{2} \sigma}\right)\right)  -\ln\left(\frac{1}{2}-\Phi\left(\frac{-U}{\sqrt{2} \sigma}\right)\right) \right].$$

Now we consider  $W(t)=\vert t\vert,\;t\in \Re,$ (i.e., the absolute error loss). Under this loss function  the usual estimator of $\theta_{(2)}$ is $\delta_{c_0}=X_{(2)}$. In this case, using Theorem 3.3, the estimator $\delta_{c_0}$ is dominated by the estimator
\begin{align*}
	\delta_{ST}(\bold{X})&=X_{(2)}-\min\bigg\{0,\frac{U}{2}\bigg\}=X_{(2)}.
\end{align*}
Also, using Theorem 3.6, the estimator $\delta_{b_0}(\bold{X})=X_{(2)}-b_0$, where $b_0$ is the solution of the equation $\int_{-\infty}^{C} \Phi\left(\frac{s}{\sigma}\right) \phi\left(\frac{s}{\sigma}\right)ds=\frac{\sigma}{4}$, is dominated by the estimator $\delta_{BZ}(\bold{X})=X_{(2)}-C,$
where $C$ is the unique solution of the following equation
\begin{equation}
	\int_{-\infty}^{C}\left[ \Phi\left(\frac{s}{\sigma}\right)- \Phi\left(\frac{s-U}{\sigma}\right)\right] \phi\left(\frac{s}{\sigma}\right)ds=\frac{\sigma}{4}-\frac{\sigma}{2}\Phi\left(\frac{-U}{\sqrt{2} \sigma}\right).
\end{equation}
\vspace*{3mm}

\noindent Now, we will illustrate an application of Theorem \ref{th1} and Corollary \ref{pitman corr}. Consider the estimation of parameter $\theta_{(2)}$ under the GPN criterion with the general loss function ($W(\cdot)$ satisfy (C1), (C2) and (C3))
\begin{eqnarray}\label{eq:3.1}
	L(\boldsymbol{\theta},a)=W(a-\theta_{(2)}), \; \boldsymbol{\theta}\in \Re^2,\; a\in \Re.
\end{eqnarray}
Under the GPN criterion, the Pitman nearest (PN) estimator of $\theta_{(2)}$ is $X_{(2)}$. In this case, using Theorems \ref{th1}, the improved estimator
$\delta_{\phi^*}(\bold{X})=X_{(2)}-\min\big\{0,\frac{U}{2}\big\}=X_{(2)}$ is same as PN estimator. Hence, there is no improvement over the estimator $X_{(2)}$ using our results.
\end{example}

\begin{example} \rm
	Let $X_1$ and $X_2$ be independent exponential random variables with $X_i$ having the p.d.f. $f(x-\theta_i),\;i=1,2,$ where
$$f(z)=\begin{cases} \frac{1}{\sigma}\, e^{-\frac{z}{\sigma}}\, ,&\text{if }\; z\geq 0\\ 0,&\text{if }\; z<0 \end{cases},$$
$\sigma>0$ is known positive constant and $\boldsymbol{\theta} \in \Re^2$ is the vector of unknown location parameters. Consider estimation of $\theta_{(2)}$ under the loss function
\begin{eqnarray}\label{eq:4.4}
	L(\boldsymbol{\theta},a)=W(a-\theta_{(2)}), \; \boldsymbol{\theta}\in \Re^2,\; a\in \Re.
\end{eqnarray}
Under the squared error loss (i.e., $W(t)=t^2,\;t\in \Re$), the usual estimator of $\theta_{(2)}$ is $X_{(2)}-\sigma$. In this case, using Theorem 3.1, the BLEE $X_{(2)}-\sigma$ is dominated by the estimator $\delta_{ST}(\bold{X})$, where
\begin{align*}
	\delta_{ST}(\bold{X})&=X_{(2)}-\min\bigg\{\sigma,\frac{\sigma}{2}+U\bigg\}=\max\bigg\{X_{(2)}-\sigma,X_{(1)}-\frac{\sigma}{2}\bigg\}.\end{align*}
 Also, using Theorem 3.4, the estimator $\delta_{b_0}(\bold{X})=X_{(2)}-\frac{3\sigma}{2}$ is dominated by the estimator $\delta_{BZ}(\bold{X})$, where
\begin{align*}
	\delta_{BZ}(\bold{X})&=X_{(2)}-\phi_{BZ}(U)=X_{(2)}-\frac{3\sigma-(2U+3\sigma)e^{-\frac{U}{\sigma}}}{2(1-e^{-\frac{U}{\sigma}})},
\end{align*}
and $U=X_{(2)}-X_{(1)}$. 


\vspace*{2mm}

Under the Linex loss function  $W(t)=e^{at}-at-1,\;t\in \Re,\; a\neq 0$, the usual estimator of $\theta_{(2)}$ is $\delta_{c_0}=X_{(2)}+\frac{1}{a} \ln(1-a\sigma)$, whenever $ a\sigma<1$. In this case, using Theorem 3.2, the estimator $X_{(2)}$ is dominated by the estimator $\delta_{ST}(\bold{X})$, where
\begin{align*}
	\delta_{ST}(\bold{X})&=X_{(2)}-\min\bigg\{-\frac{1}{a} \ln(1-a\sigma),U+\frac{\ln(2)}{a}-\frac{1}{a} \ln(2-a\sigma)\bigg\}\\
	&=\max\bigg\{X_{(2)}+\frac{1}{a} \ln(1-a\sigma),X_{(1)}+\frac{1}{a} \ln(2-a\sigma)-\frac{\ln(2)}{a}\bigg\}.
\end{align*}
	Also, using Theorem 3.5 the estimator $\delta_{b_0}(\bold{X})=X_{(2)}-\frac{1}{a}[\ln2 -\ln(1-a\sigma)-\ln(2-a\sigma)]$, whenever $ a\sigma<1$, is dominated by the estimator $\delta_{BZ}(\bold{X})$, where
	$$ \delta_{BZ}(\bold{X})=X_{(2)}-\phi_{BZ}(U)=X_{(2)}-\frac{1}{a}\left[\ln2 -\ln(1-a\sigma)-\ln(2-a\sigma)+\ln\left(1-e^{-U\left(\frac{1}{\sigma}-a\right)}\right)-\ln\left(1-e^{-\frac{U}{\sigma}}\right)\right].$$

Under the absolute error loss $W(t)=\vert t\vert,\;t\in \Re,$ the usual estimator of $\theta_{(2)}$ is $\delta_{c_0}=X_{(2)}-\sigma \ln(2)$.  In this case, using Theorem 3.3, the estimator $X_{(2)}-\sigma \ln(2)$ is dominated by the estimator
\begin{align*}
	\delta_{ST}(\bold{X})&=X_{(2)}-\min\bigg\{\sigma \ln(2),U+\frac{\sigma}{2}\ln(2)\bigg\}=\max\bigg\{X_{(2)}-\sigma \ln(2),X_{(1)}-\frac{\sigma}{2}\ln(2)\bigg\}.
\end{align*}
Also, using Theorem 3.6, the estimator $\delta_{b_0}(\bold{X})=X_{(2)}+\sigma\ln\left(1-\frac{1}{\sqrt{2}}\right)$ is dominated by the estimator $\delta_{BZ}(\bold{X})=X_{(2)}-C,$
where $C>0$ is the unique solution of the following equation
\begin{align*}
\int_{0}^{C}\int_{0}^{\min\{U,s\}} e^{-\frac{2s}{\sigma}}e^{\frac{y}{\sigma}}dy\,ds=\frac{\sigma^2}{4}\left(1-e^{-\frac{U}{\sigma}}\right).
\end{align*}


\vspace*{3mm}

\noindent Now, we will illustrate an application of Theorem \ref{th1} and Corollary \ref{pitman corr}. Consider the estimation of parameter $\theta_{(2)}$ under the GPN criterion with the general loss function ($W(\cdot)$ satisfy (C1), (C2) and (C3))
\begin{eqnarray}\label{eq:3.1}
	L(\boldsymbol{\theta},a)=W(a-\theta_{(2)}), \; \boldsymbol{\theta}\in \Re^2,\; a\in \Re.
\end{eqnarray}
Under the GPN criterion, the PN estimator of $\theta_{(2)}$ is $X_{(2)}-\sigma \ln(2)$. In this case, using Theorems \ref{th1}, the estimator
$$\delta_{\phi_0}(\bold{X})=\max\bigg\{X_{(2)}-\sigma \ln(2),X_{(1)}-\frac{\sigma}{2}\ln(2)\bigg\}$$ is the Pitman nearer to $\theta_{(2)}$ than the estimator $X_{(2)}-\sigma \ln(2)$, under the GPN criterion.

\end{example}

\section{Simulation study}\label{sec7}
\begin{figure}[h!]
	\begin{subfigure}{.5\textwidth}
		\centering
		\includegraphics[width=80mm,scale=1.2]{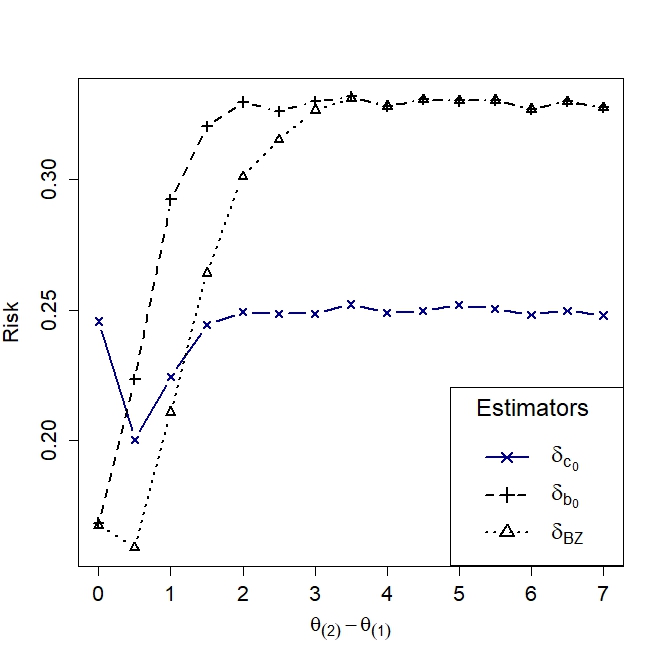} 
		\caption{$\sigma=0.5$} 
		
	\end{subfigure}
	\begin{subfigure}{.5\textwidth}
		\centering
		\includegraphics[width=80mm,scale=1.2]{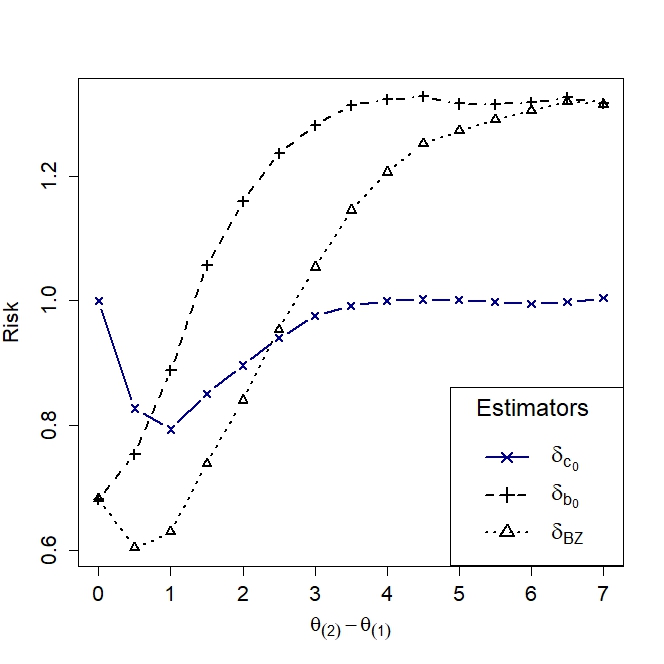} 
		
		\caption{$\sigma=1$} 
		
	\end{subfigure}
	\\	\begin{subfigure}{.5\textwidth}
		\centering
		\includegraphics[width=80mm,scale=1.2]{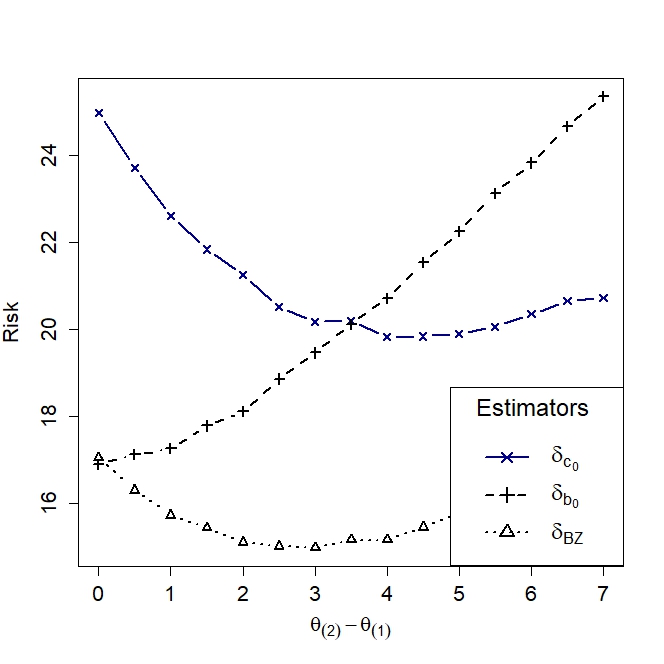} 
		\caption{$\sigma=5$} 
		
	\end{subfigure}
	\begin{subfigure}{.5\textwidth}
		\centering
		
		\includegraphics[width=80mm,scale=1.2]{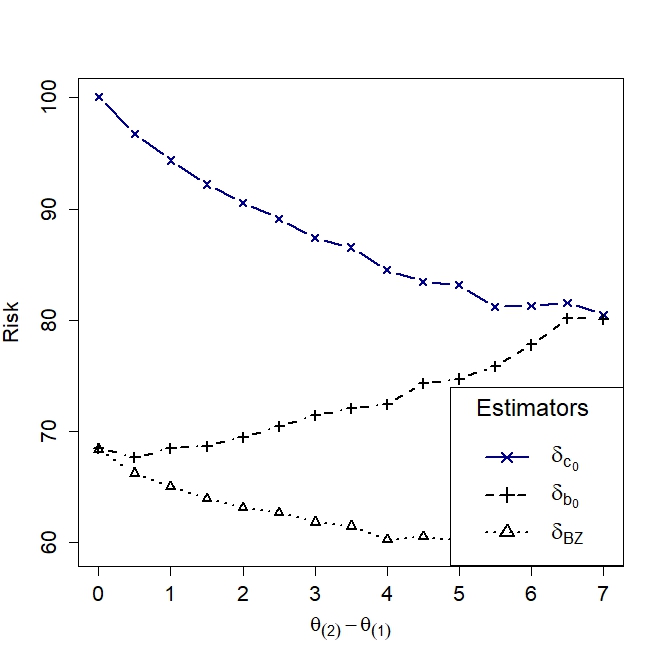} 
		
		\caption{$\sigma=10$} 
		
	\end{subfigure}

	\caption{Risk plots of estimators of parameter $\theta_{(2)}$ against values of $\theta=\theta_{(2)}-\theta_{(1)}$: when $W(t)=t^2,\;t\in \Re,$ and so, $b_0=\frac{\sigma}{\sqrt{\pi}}$.}
	\label{fig1}
\end{figure}

\begin{figure}
	\begin{subfigure}{.5\textwidth}
		\centering
		\includegraphics[width=70mm,scale=1]{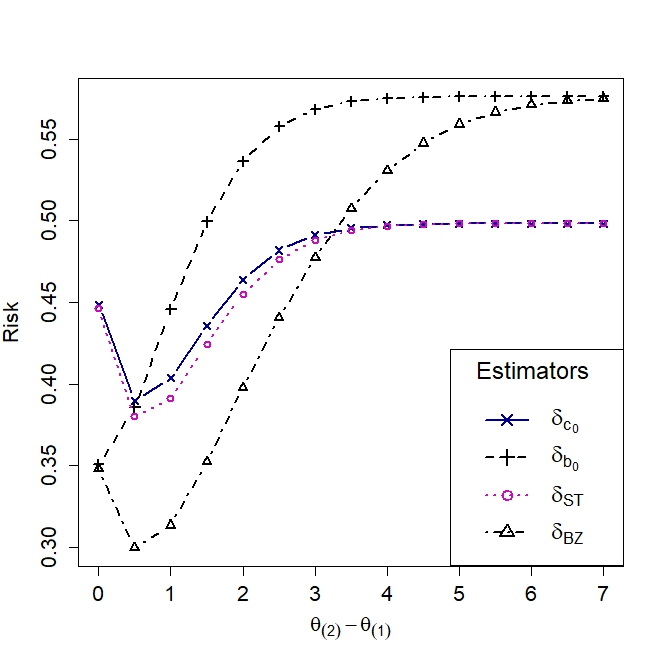} 
		\caption{$\sigma=1$ and $a=1$} 
		
	\end{subfigure}
	\begin{subfigure}{.5\textwidth}
		\centering
		\includegraphics[width=70mm,scale=1]{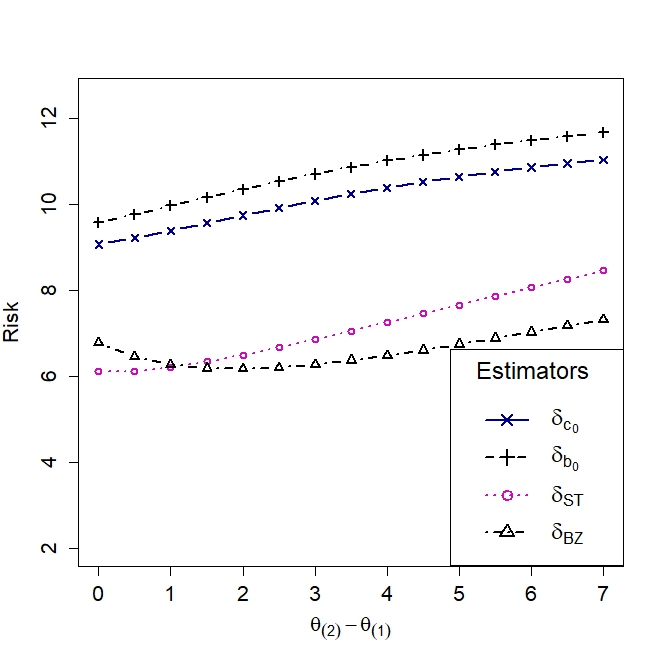} 
	
		\caption{$\sigma=5$ and $a=1$} 
		
	\end{subfigure}
	\\	\begin{subfigure}{.5\textwidth}
		\centering
		\includegraphics[width=70mm,scale=1.2]{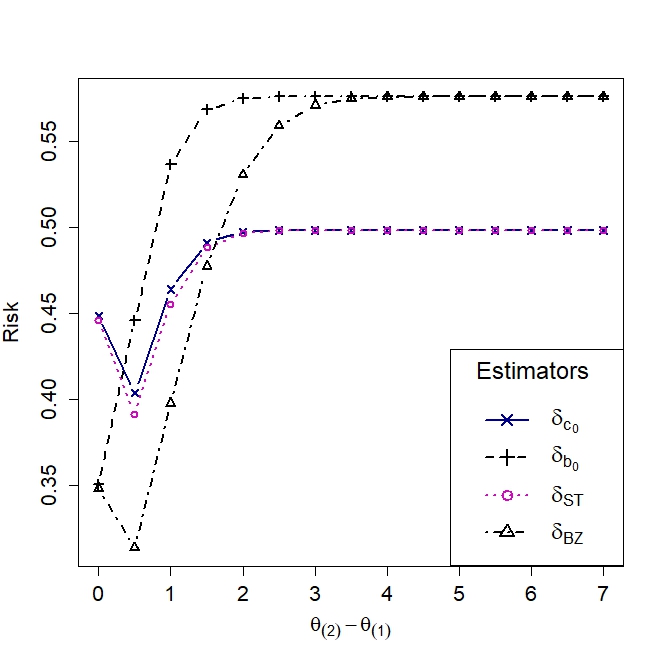} 
		\caption{$\sigma=0.5$ and $a=2$} 
		
	\end{subfigure}
	\begin{subfigure}{.5\textwidth}
		\centering
		
		\includegraphics[width=70mm,scale=1.2]{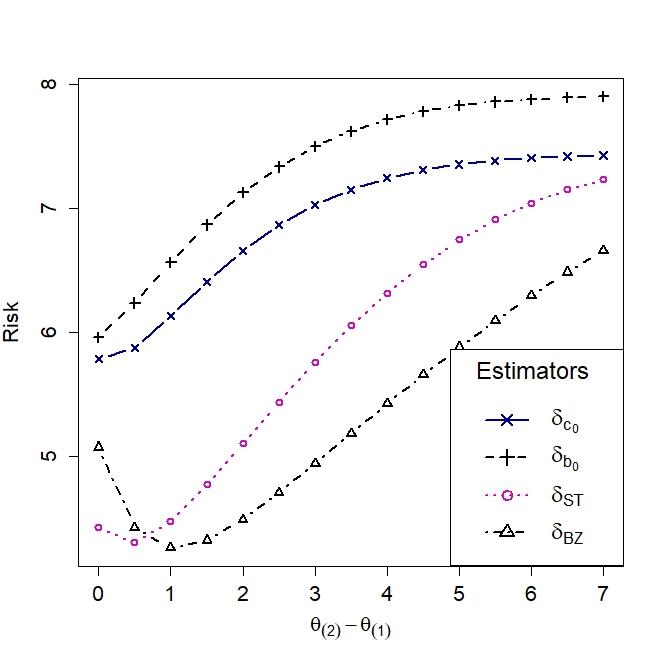} 
		
		\caption{$\sigma=2$ and $a=2$} 
		
	\end{subfigure}
	\\	\begin{subfigure}{.5\textwidth}
	\centering
	\includegraphics[width=70mm,scale=1.2]{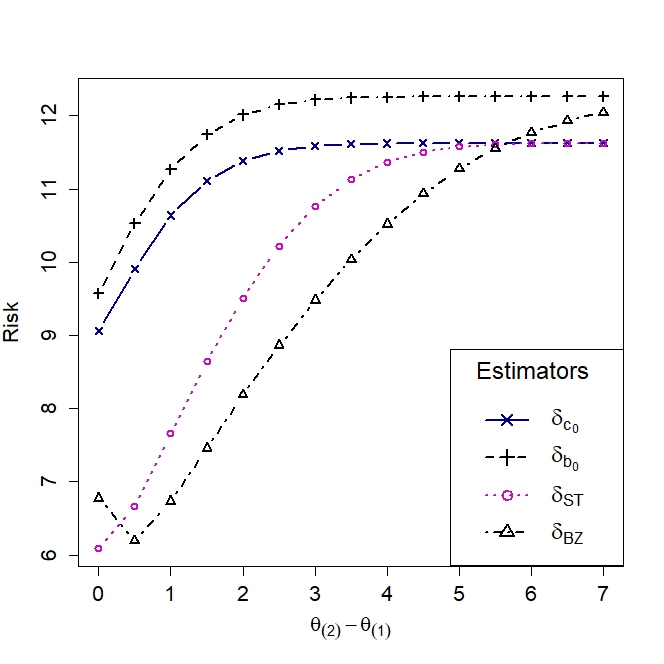} 
	\caption{$\sigma=1$ and $a=5$} 
	
\end{subfigure}
\begin{subfigure}{.5\textwidth}
	\centering
	
	\includegraphics[width=70mm,scale=1.2]{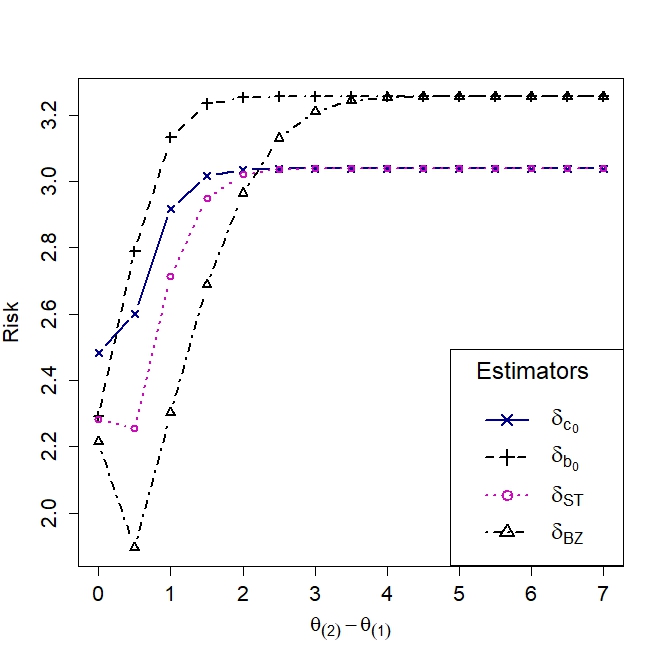} 
	
	\caption{$\sigma=0.5$ and $a=5$} 
	
\end{subfigure}
	
	\caption{Risk plots of estimators of parameter $\theta_{(2)}$ against values of $\theta=\theta_{(2)}-\theta_{(1)}$: when $W(t)=e^{at}-at-1,\;t\in \Re,\; a\neq 0$, and so, $b_0=\frac{1}{a}\left[\ln2+\frac{a^2\sigma^2}{2} +\ln\left(\Phi\left(\frac{a\sigma}{\sqrt{2}}\right)\right)\right]$.}
	\label{fig2}
\end{figure}

\begin{figure}
	\begin{subfigure}{.5\textwidth}
		\centering
		\includegraphics[width=80mm,scale=1.2]{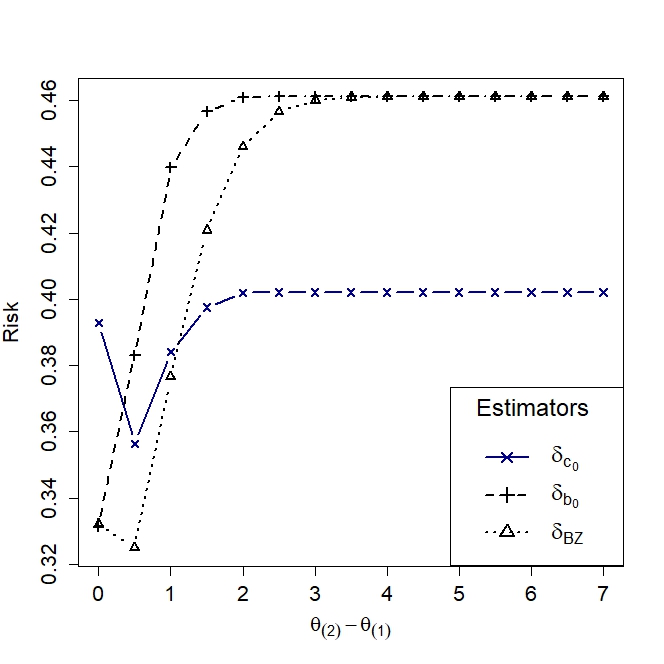} 
		\caption{$\sigma=0.5$} 
		
	\end{subfigure}
	\begin{subfigure}{.5\textwidth}
		\centering
		\includegraphics[width=80mm,scale=1.2]{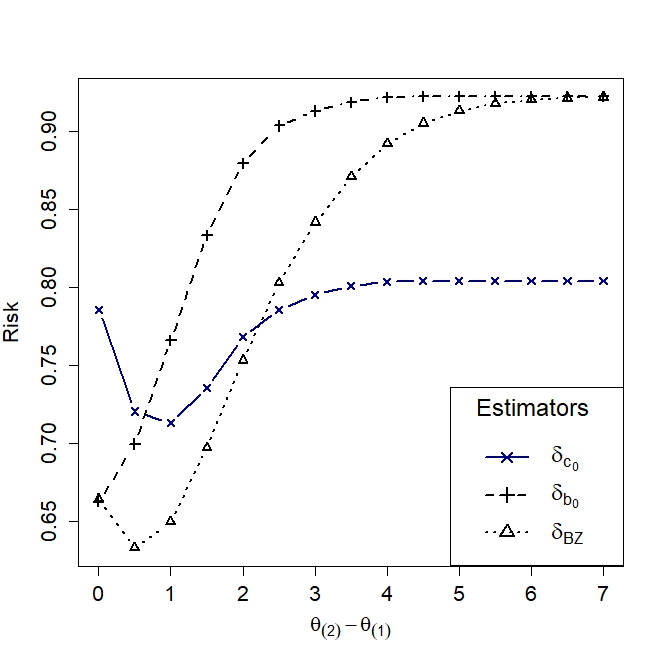} 
		
		\caption{$\sigma=1$} 
		
	\end{subfigure}
	\\	\begin{subfigure}{.5\textwidth}
		\centering
		\includegraphics[width=80mm,scale=1.2]{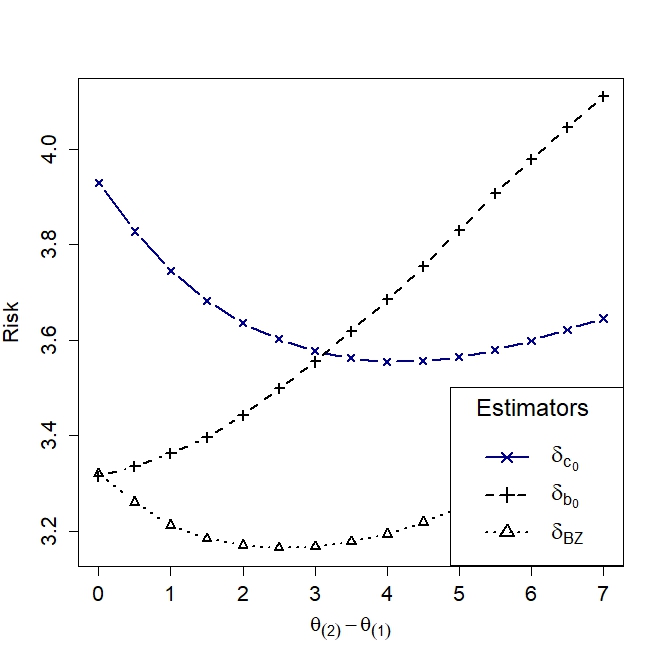} 
		\caption{$\sigma=5$} 
		
	\end{subfigure}
	\begin{subfigure}{.5\textwidth}
		\centering
		
		\includegraphics[width=80mm,scale=1.2]{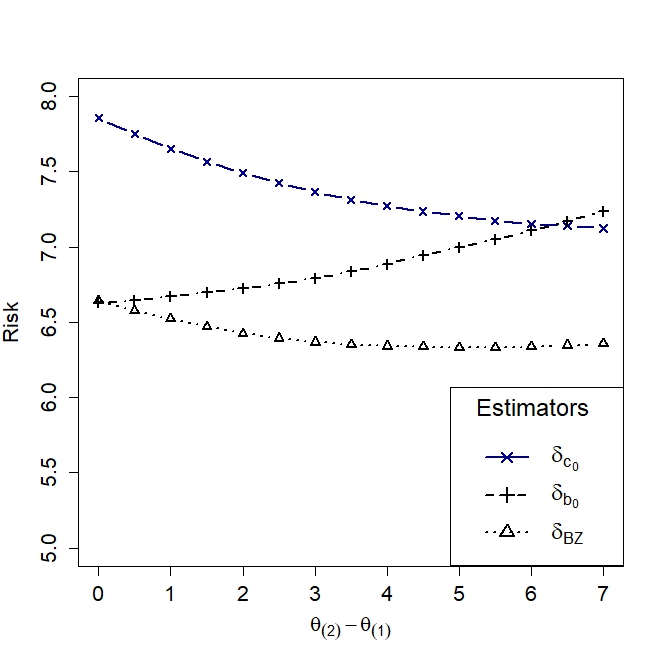} 
		
		\caption{$\sigma=10$} 
		
	\end{subfigure}
	
	\caption{Risk plots of estimators of parameter $\theta_{(2)}$ against values of $\theta=\theta_{(2)}-\theta_{(1)}$: when $W(t)=\vert t\vert,\;t\in \Re,$ and so, $b_0$ is the solution of the equation $\int_{-\infty}^{C} \Phi\left(\frac{s}{\sigma}\right) \phi\left(\frac{s}{\sigma}\right)ds=\frac{\sigma}{4}$.}
	\label{fig3}
\end{figure}

In Example 4.1, we considered the estimation of the larger mean $\theta_{(2)}$ of two independent normal distributions with a known common variance ($\sigma^2$). We derived estimators that improved upon the usual estimator $\delta_{c_0}(\boldsymbol{X}) = X_{(2)}$ and the estimator $\delta_{1}(\boldsymbol{X}) = X_{(2)} - b_0$ under three different loss functions. To evaluate the risk performance of these estimators of $\theta_{(2)}$ under these loss functions, we conducted Monte Carlo simulations to compare the risk performances of the usual estimator $X_{(2)}$, the Stein-type estimator $\delta_{ST}$, the estimator $\delta_{1}(\boldsymbol{X})$, and the Brewster-Zidek type estimator $\delta_{BZ}$. We computed the simulated risks based on 50000 simulations from relevant distributions, and the resulting risk function of the different estimators were plotted in Figures \ref{fig1}-\ref{fig3}. The following observations are evident from Figures \ref{fig1}-\ref{fig3}:

\begin{itemize}
\item[(i)] The risk of the Stein type estimator $\delta_{ST}$ has smaller risk than the usual estimator $\delta_{c_0}$. For smaller $\sigma$, it can be seen that the risk values of both the estimators are almost equal when $\theta_{(2)}-\theta_{(1)}$ takes values close to zero and approximately more than 3. For higher values of $\sigma$, we observe that $\delta_{ST}$ performs significantly better than $\delta_{c_0}$. In this case, the region on improvement is better than the previous (small $\sigma$) case. 
\item[(ii)]  The risk performance of $\delta_{BZ}$ is better then $\delta_{1}$ for all $\theta_{(2)}-\theta_{(1)}$. It is observed that for larger values of $\sigma$, the improvement interval is bigger compared to the smaller values of $\sigma$. 

\item[(iii)] We also found that there was no clear winner among the various estimators, as the estimator $\delta_{1}$ and the $\delta_{BZ}$ performed better than the usual estimator and $\delta_{ST}$ for small and moderate values of $\theta_{(2)}-\theta_{(1)}$, whereas the estimator $\delta_{ST}$ dominated the other two estimators for large values of $\theta_{(2)}-\theta_{(1)}$.
\end{itemize}

\FloatBarrier
\FloatBarrier
\FloatBarrier
\FloatBarrier

\section{\textbf{Real-life Data Analysis}}
\label{sec8}
\noindent

\vspace*{3mm}

For real life data analysis, we have considered the ``Jute fiber breaking strength data", discussed by \cite{xia2009study} and presented in Table \ref{table:1}. This data represents the breaking strength of jute fibre of two different gauge lengths. 

\vspace*{2mm}

Jute is a versatile natural fiber widely employed in various products, including textiles, ropes, sacks, and geotextiles. The breaking strength of jute fibers plays a pivotal role in determining the quality and durability of these products. By estimating the maximum breaking strength with different gauge lengths, manufacturers can ensure that their products conform to the required standards and can withstand a range of stresses and loads. Therefore, the estimation of the maximum breaking strength of jute fibers using two different gauge lengths holds significant importance in the realm of material science and engineering for several compelling reasons.

\vspace*{2mm}

For data analysis, we initially examined whether these datasets follow two-parameter exponential distributions. We used the one-sample Kolmogorov-Smirnov test to see if the data with gauge lengths 10 mm and 15 mm may be following a two parameter exponential distributions and found p-values of $0.755$ and $0.306$, respectively, indicating that the data with gauge lengths 10 mm and 15 mm follow a two parameter exponential distributions with a common scale parameter value of $322$.

\begin{table}[h!]
	\caption{The breaking strength of jute fibre}
	\centering
	\begin{tabular}{@{}|c| c| c| @{}}
		\hline	No. & gauge length 10 mm & gauge length 15 mm \\
		 \hline  \hline
		1   & 43.93              & 594.4              \\
		2   & 50.16              & 202.75             \\
		3   & 101.15             & 168.37             \\
		4   & 108.94             & 574.86             \\
		5   & 123.06             & 225.65             \\
		6   & 141.38             & 76.38              \\
		7   & 151.48             & 156.67             \\
		8   & 163.4              & 127.81             \\
		9   & 177.25             & 813.87             \\
		10  & 183.16             & 562.39             \\
		11  & 212.13             & 468.47             \\
		12  & 257.44             & 135.09             \\
		13  & 262.9              & 72.24              \\
		14  & 291.27             & 497.94             \\
		15  & 303.9              & 355.56             \\
		16  & 323.83             & 569.07             \\
		17  & 353.24             & 640.48             \\
		18  & 376.42             & 200.76             \\
		19  & 383.43             & 550.42             \\
		20  & 422.11             & 748.75             \\
		21  & 506.6              & 489.66             \\
		22  & 530.55             & 678.06             \\
		23  & 590.48             & 457.71             \\
		24  & 637.66             & 106.73             \\
		25  & 671.49             & 716.3              \\
		26  & 693.73             & 42.66              \\
		27  & 700.74             & 80.4               \\
		28  & 704.66             & 339.22             \\
		29  & 727.23             & 70.09              \\
		30  & 778.17             & 193.42            \\
		\hline
	\end{tabular}
	\label{table:1}
\end{table}

\vspace*{2mm}

Let $X_1$ and $X_2$ be the two independent random variables representing the breaking strength of jute fibre corresponding to gauge length 10 mm and 15 mm, respectively. Therefore $X_1\sim Exp(\theta_1,\sigma)$ and $X_2\sim Exp(\theta_2,\sigma)$, where $\theta_1$ and $\theta_2$ are location parameters, and common known scale parameter $\sigma= 322/30 = 10.73$. Using our finding of this paper, we obtain estimates that are better than the natural estimates of parameter $\max\{\theta_1,\theta_2\}$.
\vspace*{2mm}

In Example 5.2, we have presented various estimators under the squared error loss, linex loss and absolute error loss functions, which are given in Table \ref{table:2}, Table \ref{table:3} and Table \ref{table:4}, respectively. From theoretical results (in Example 5.2), we infer that the Stein type estimated value $\delta_{ST}(\bold{x})$ is better than estimated value of $\delta_{c_0}(\bold{x}) $ and the Brewster-Zidek type estimated value $\delta_{BZ}(\bold{x})$ is better than estimated value of $\delta_{b_0}(\bold{x}) $ for parameter $\max\{\theta_1,\theta_2\}$.

\FloatBarrier
\begin{table}[h!]
	\caption{Various estimated values of parameter $\max\{\theta_1,\theta_2\}$: \\under the squared error loss}
	\centering
	\begin{tabular}{@{}|c |c |c |c| @{}}
		\hline
		$\;\;\delta_{c_0}(\bold{x})\;\;$ & $\quad\delta_{ST}(\bold{x}) \quad$ & $\quad\delta_{b_0}(\bold{x}) \quad$ & $\quad \delta_{BZ}(\bold{x}) \quad$
		\\ \hline \hline
	33.2   & 37.3   & 27.835 & 37.94                                          \\ \hline
		
	\end{tabular}
	\label{table:2}
\end{table}
\FloatBarrier

\FloatBarrier
\begin{table}[h!]
	\caption{Various estimated values of parameter $\max\{\theta_1,\theta_2\}$:\\ under the Linex loss with $a=-1$}
	\centering
	\begin{tabular}{@{}|c |c |c |c| @{}}
		\hline
		$\;\;\delta_{c_0}(\bold{x})\;\;$ & $\quad\delta_{ST}(\bold{x}) \quad$ & $\quad\delta_{b_0}(\bold{x}) \quad$ & $\quad \delta_{BZ}(\bold{x}) \quad$
		\\ \hline \hline
		41.47   & 41.47   & 39.62 & 41.52                                          \\ \hline
		
	\end{tabular}
	\label{table:3}
\end{table}
\FloatBarrier

\FloatBarrier
\begin{table}[h!]
	\caption{Various estimated values of parameter $\max\{\theta_1,\theta_2\}$: \\under the absolute error loss}
	\centering
	\begin{tabular}{@{}|c |c |c |c| @{}}
	\hline
	$\;\;\delta_{c_0}(\bold{x})\;\;$ & $\quad\delta_{ST}(\bold{x}) \quad$ & $\quad\delta_{b_0}(\bold{x}) \quad$ & $\quad \delta_{BZ}(\bold{x}) \quad$
	\\ \hline \hline
		36.5   & 38.94    & 30.75 & 45.65                                            \\ \hline
		
	\end{tabular}
	\label{table:4}
\end{table}
\FloatBarrier

\section{Conclusions}\label{sec9}
In the present article, we have considered the estimation, from a decision-theoretic point of view, of the larger location parameter of two general location family of distributions under a general location invariant loss function. We have proposed a Stein-type estimator, which improves upon the usual estimator. Next, we have proposed an estimator, $\delta_{b_0}$, which is similar to the usual estimator. Using the IERD approach of Kubokawa, a class of estimators has been derived that dominates $\delta_{b_0}$. It is seen that the boundary estimator of this class is the Brewster-Zidek type estimator. As an application, the improved estimators are derived for particular loss functions. We have also considered the same estimation problem with respect to the generalized Pitman nearness criterion. We have proposed an estimator that is nearer than the usual estimator. As an application, explicit expressions of all the improved estimators are obtained for the normal and exponential models. Further, we have compared the risk performance of the proposed estimators using the Monte Carlo simulation. Finally, we present a real-life application that demonstrates the practical significance of the findings presented in this paper.
\begin{small}
	\bibliography{references}		
\end{small}
\end{document}